\documentclass[10pt]{amsart}
\usepackage{latexsym, amsmath}

\newcommand{\varespilon}{\varepsilon}

\renewcommand{\epsilon}{\varepsilon}
\newcommand{\newsection}[1]
{\subsection{#1}\setcounter{theorem}{0} \setcounter{equation}{0}
\par\noindent}

\newtheorem{theorem}{Theorem}

\newtheorem{lemma}[theorem]{Lemma}
\newtheorem{corr}[theorem]{Corollary}

\newtheorem{proposition}[theorem]{Proposition}
\newtheorem{deff}[theorem]{Definition}

\newcommand{\bth}{\begin{theorem}}
\newcommand{\ble}{\begin{lemma}}
\newcommand{\bcor}{\begin{corr}}

\newcommand{\bdeff}{\begin{deff}}

\newcommand{\bprop}{\begin{proposition}}
\newcommand{\eth}{\end{theorem}}
\newcommand{\ele}{\end{lemma}}
\newcommand{\ecor}{\end{corr}}
\newcommand{\edeff}{\end{deff}}

\newcommand{\eprop}{\end{proposition}}

\newcommand{\cd}{\, \cdot\, }

\renewcommand{\Pi}{\varPi}

\renewcommand{\epsilon}{\varepsilon}

\newcommand{\parital}{\partial}

\newcommand{\R}{{\mathbb R}}

\begin{document}

\title[Global existence for wave equations]
{Hyperbolic trapped rays and global existence of quasilinear wave
equations}
\thanks{The authors were supported in part by the NSF}

\author{Jason Metcalfe}
\address{School of Mathematics, Georgia Institute of Technology,
Atlanta, GA  30332-0160}
\author{Christopher D. Sogge}
\address{Department of Mathematics,  Johns Hopkins University,
Baltimore, MD 21218}

\maketitle

\newsection{Introduction}

The purpose of this paper is to give a simple proof of global
existence for quadratic quasilinear Dirichlet-wave equations
outside of a wide class of compact obstacles in the critical case
where the spatial dimension is three.  Our results improve on
earlier ones in Keel, Smith and Sogge \cite{KSS} in several ways.
First, and most important, we can drop the star-shaped hypothesis
and handle non-trapping obstacles as well as any obstacle that has
exponential local decay rate of energy for $H^2$ data for the
linear equation (see \eqref{energydecay} below). This hypothesis
is fulfilled in the non-trapping case where there is actually
exponential local decay of energy \cite{MRS} with no loss of
derivatives. This hypothesis \eqref{energydecay} is also known to
hold in several examples involving hyperbolic trapped rays. For
instance, our results apply to  situations where the obstacle is a
finite union of convex bodies with smooth boundary (see
\cite{Ikawa1}, \cite{Ikawa2}). In addition to improving the
hypotheses on the obstacles, we can also improve considerably on
the decay assumptions on the initial data at infinity compared to
the results in \cite{KSS} which were obtained by the conformal
method. Lastly, we are able handle non-diagonal systems involving
multiple wave speeds.

We shall use a refinement of techniques developed in earlier work
of Keel, Smith and Sogge \cite{KSS2}, \cite{KSS3}.  In particular,
we shall use a modification of Klainerman's commuting vector
fields method \cite{knull} that only uses the collection of vector
fields that seems ``admissible" for boundary value problems.

The main innovation in this approach versus the classical one for
the boundaryless case is the use of weighted space-time $L^2$
estimates to handle the various lower order terms that necessarily
arise in obstacle problems. The weights involved are just negative
powers of $\langle x\rangle$. These couple well with the pointwise
estimates that we use, which involve $O(\langle x\rangle^{-1})$
decay of solutions of linear inhomogeneous Dirichlet-wave
equations, as opposed to the more standard $O(t^{-1})$ decay for
the boundaryless case, which are much more difficult to obtain for
obstacle problems. Because of the fact that we are dealing with
such problems, it does not seem that we can use vector fields such
as the generators of hyperbolic rotations,
$x_i\partial_t+t\partial_i$, $i=1,2,3$. Additionally, it seems
that these cannot be used for multiple wave speed problems since
they have an associated speed (one in the above case).  So, unlike
in Klainerman's argument \cite{knull} for the Minkowski space
case, we are only able to use the generators of spatial rotations
and space-time translations
\begin{equation}\label{1.1}
Z=\{ \, \partial_i, x_j\partial_k-x_k\partial_j, \, 0\le i\le 3,
1\le j<k\le 3\},\end{equation} as well as the scaling vector field
\begin{equation}\label{1.2}
L=t\partial_t+r\partial_r.
\end{equation}
Here, and in what follows, we are using the notation that
$(x_1,x_2,x_3)$ denote the spacial coordinates, while either $x_0$
or $t$ will denote the time coordinate, depending on the context.
Also, $r=|x|$, and $\langle x\rangle = \langle r\rangle =
\sqrt{1+r^2}$.  We shall also let $\partial = \partial_{t,x}$
denote the space-time gradient.

Another difficulty that we encounter in the obstacle case is
related to the simple fact that while the vector fields
\begin{equation}\label{1.3}
\Omega_{ij}=x_i\partial_j-x_j\partial_i, \quad 1\le i<j\le 3
\end{equation}
and $L$ preserve the equation $(\partial_t^2-\Delta)u=0$ in the
Minkowski space case if $u$ is replaced by either $Lu$ or
$\Omega_{ij}u$, this is not true in the obstacle case due to the
fact that the Dirichlet boundary conditions are not preserved by
these operators.  Since the generators of spatial rotations,
$\Omega_{ij}$, have coefficients that are small near our compact
obstacle, this fact is somewhat easy to get around when dealing
with them; however, it is a bit harder to deal with the scaling
vector field, $L$, since its coefficients become large on the
obstacle as $t$ goes to infinity. As a result, we are forced to
consider in our estimates combinations of the $Z$ operators and
the $L$ operators that involve relatively few of the scaling
vector fields.  This, together with the fact that there is
necessarily a loss of smoothness in the local energy estimates for
obstacles with trapped rays, makes the combinatorics that arise
more complicated than in the Minkowski space case first studied by
Klainerman \cite{knull}.

In earlier works \cite{KSS}, \cite{KSS3} the obstacle was assumed
to be star-shaped.  This was a convenient assumption in proving
energy estimates involving the scaling operator $L$.  For
instance, in proving energy estimates for $Lu$ for solutions of
$(\partial^2_t-\Delta)u=0$ one finds that if $\mathcal{K}$ is
star-shaped then, although energy is not conserved, the
contribution from the boundary to energy identities has a
favorable sign.  This is in the spirit of Morawetz's original
argument \cite{M}.  If one drops the star-shaped assumption this
argument of course breaks down.  However, in this paper we exploit
the fact that we still can prove favorable estimates for solutions
of {\it nonlinear equations}.  The additional terms arising from
the boundary can be estimated using Lemma \ref{lemma2.9} and the
$L^2_tL^2_x(\langle x\rangle^{-1/2}dxdt)$ estimates since the
forcing terms are nonlinear functions of $(du,du^2)$ that vanish
to second order.

Let us now describe more precisely our assumptions on our
obstacles $\mathcal{K}\subset \mathbb{R}^3$.  We shall assume that
$\mathcal{K}$ that is smooth and compact. We do not assume that
$\mathcal{K}$ is connected. Without loss of generality, we may
assume throughout that
$$\mathcal{K}\subset \{x\in \mathbb{R}^3: \, |x|<1\}.$$
Our only additional assumption is that there is exponential local
decay of energy with a possible loss of derivatives.  To be
specific, we require that there be a $c>0$, a constant $C$  so
that
\begin{multline}\label{energydecay}
\Bigl(\int_{\{x\in \mathbb{R}^3\backslash \mathcal{K}: \,
|x|<4\}}|u'(t,x)|^2\, dx\Bigr)^{1/2} \le Ce^{-c t}
\sum_{|\alpha|\le 1} \|\partial_x^\alpha u'(0,\cd)\|_2, \\
\text{if } \, u(0,x)=\partial_tu(0,x)=0, \, \, \text{for } \,
|x|>4.
\end{multline}

We remark that our results do not actually require exponential
decay of local energy.  A decay rate of $O(\langle
t\rangle^{-3-\delta})$, $\delta>0$ would suffice since our main
$L^2$-estimates involve 3 or fewer powers of the scaling operator
$L$.    By tightening the arguments one might even be able to show
that $O(\langle t\rangle^{-1-\delta})$ is sufficient.  On the
other hand, we shall assume \eqref{energydecay} throughout since
the proofs under this weaker decay rate would be more technical.
Moreover, if in the 3-dimensional case, all of the examples that
we know that involve polynomial decay actually have exponential
decay of local energy. For related problems in general relativity,
though, it might be much easier to establish local polynomial
decay of energy.

Recall that if the obstacle is star-shaped or non-trapping a
stronger version of  \eqref{energydecay} is always valid where in
the right one just takes the $H^1\times L^2$-norm of
$(u(0,\cd),\partial_tu(0,\cd))$ (see Lax, Morawetz and Phillips
\cite{LMP} for the star-shaped case, and Morawetz, Ralston and
Strauss \cite{MRS} for the non-trapping case, and Melrose
\cite{Mel} for further results of this type). On the other hand,
if $\mathbb{R}^3\backslash \mathcal{K}$ contains any trapped rays,
then Ralston \cite{ralston} showed that this stronger inequality
cannot hold. So there must be some {\it ``loss"} $\ell>0$ of
regularity if there is energy decay when there are trapped rays.
In \eqref{energydecay} we are assuming that $\ell=1$.  By
interpolation, there is no loss of generality in making this
assumption since if the analog of \eqref{energydecay} held where
the sum was taken over a given $\ell >1$ then \eqref{energydecay}
would still be valid (with a different constant in the
exponential).  (The same argument shows that a variant of
\eqref{energydecay} holds where one has
$\|u(0,\cd)\|_{H^{1+\delta}_D}+\|\partial_t
u(0,\cd)\|_{H^{\delta}_D}$ in the right for $\delta>0$.)

 In other direction, Ikawa \cite{Ikawa1},
\cite{Ikawa2} was able to show that if $\mathcal{K}$ is a finite
union of convex obstacles with smooth boundary then one has
exponential decay of local energy with a loss of $\ell=7$
derivatives, which as we just pointed out leads to
\eqref{energydecay} here.  Ikawa's theorem requires additional
technical assumptions that we shall not describe (see
\cite{Ikawa2}); however, they are always satisfied for instance in
the case where $\mathcal{K}$ is the union of two disjoint convex
obstacles or any number of balls that are sufficiently separated.
Thus, even for the case where $\mathcal{K}$ is the union of $3$
sufficiently separated balls one can always have infinitely many
trapped rays and still have \eqref{energydecay} (and the nonlinear
results to follows).  We also mention the work of  Burq
\cite{burq}  who showed that for {\it any} compact obstacle
$\mathcal{K}$ with smooth boundary, one has a local decay that is
$O((\log(2+t))^{-k})$ for any $k$ if one takes the loss of
regularity to be $\ell=k$. Such a decay rate is not fast enough
for us to be able to prove global existence for this class of
obstacles, and it seems doubtful that such results could hold in
this context since Burq's results include the case where
$\mathbb{R}^3\backslash\mathcal{K}$ has trapped elliptic rays.  On
the other hand, an interesting question would be whether our
hypothesis \eqref{energydecay} might hold under the assumption
that $\mathbb{R}^3\backslash \mathcal{K}$ only contains hyperbolic
trapped rays.

  For obstacles  ${\mathcal K}\subset
{\mathbb R}^3$, as above satisfying \eqref{energydecay} we shall
consider smooth, quadratic, quasilinear systems of the form
\begin{equation}\label{1.5}
\begin{cases}
\square u = Q(du,d^2u), \quad (t,x)\in {\mathbb R}_+\times
{\mathbb R}^3\backslash {\mathcal K}
\\
u(t,\cdot)|_{\partial\mathcal K}=0
\\
u(0,\cdot)=f, \, \, \partial_tu(0,\cdot)=g.
\end{cases}
\end{equation}
Here,
\begin{equation}\label{1.6}
\square = (\square_{c_1},\square_{c_2},\dots,\square_{c_D}),
\end{equation} is a vector-valued multiple speed D'Alembertian
with
$$\square_{c_I}=\partial^2_t-c^2_I\Delta,$$
where we assume that the wave speeds $c_I$ are all positive but
not necessarily distinct.  Also, $\Delta =
\partial_1^2+\partial_2^2+\partial_3^2$ is the standard
Laplacian.  By a simple scaling argument, in showing that
\eqref{1.5} admits global small amplitude solutions, as mentioned
before, we shall assume without loss of generality that
$\mathcal{K}\subset \{x \in \mathbb{R}^3: \, |x|<1\}$.

By quasilinear we mean that the nonlinear term $Q(du,d^2u)$ is
linear in the second derivatives of $u$.  We shall also assume
that the highest order nonlinear terms are symmetric, by which we
mean that, if we let $\partial_0=\partial_t$, then
\begin{equation}\label{1.7}
Q^I(du,d^2u)=B^I(du)+\sum_{\substack{0\le j,k,l\le 3\\ 1\le J,K\le
D} }B^{IJ,jk}_{K,l}\partial_lu^K \,
\partial_j\partial_ku^J, \quad 1\le I\le D,
\end{equation}
with $B^I(du)$ a quadratic form in the gradient of $u$, and
$B^{IJ,jk}_{K,l}$ real constants satisfying the symmetry
conditions
\begin{equation}\label{1.8}
B^{IJ,jk}_{K,l}=B^{JI,jk}_{K,l}=B^{IJ,kj}_{K,l}.
\end{equation}

To obtain global existence, we shall also require that the
equations satisfy a form of the null condition of Christodoulou
and Klainerman.  Let us first assume, for simplicity that the wave
speeds $c_I$, $I=1,\dots,D$ are distinct.  In this case, the null
condition for the quasilinear terms only involves
self-interactions of each wave family.  Specifically, we require
that self-interactions among the quasilinear terms satisfy the
standard null condition for the various wave speeds:
\begin{equation}\label{1.9}
\sum_{0\le j,k,l\le 3}B^{IJ,jk}_{J,l}\xi_j\xi_k\xi_l=0 \quad
\text{whenever}\quad
\frac{\xi_0^2}{c^2_J}-\xi_1^2-\xi_2^2-\xi_3^2=0, \quad
I,J=1,\dots,D.
\end{equation}

For the quasilinear terms, if one allows repeated wave speeds, it
will be required that the interactions of families with the same
speed satisfy a null condition.  Specifically if we let
$\mathcal{I}_p=\{I\: : \: c_I=c_{I_p},\: 1\le I\le D\}$ then the
above null condition is extended to
$$\sum_{j,k,l\le
3}B^{IJ,jk}_{K,l}\xi_j\xi_k\xi_l=0\:\text{ whenever }\:
\frac{\xi_0^2}{c^2_{I_p}}-\xi_1^2-\xi^2_2-\xi_3^2=0,\:
(J,K)\in\mathcal{I}_p\times\mathcal{I}_p,\, 1\le I\le D.$$

To describe the null condition for the lower order terms, we note
that we can expand
$$B^I(du)=\sum_{\substack{1\le J,K\le D\\ 0\le j,k\le 3}}
A^{I,jk}_{JK}\partial_ju^J\partial_ku^K.$$ We then require that
each component satisfy the standard null condition for multiple
wave speeds
\begin{equation}\label{1.10}
\sum_{0\le j,k\le 3} A^{I,jk}_{JK}\xi_j\xi_k=0\quad \text{for all
}\quad \xi \in \R \times \R^3,  \quad 1\le J,K\le D.
\end{equation}
This means that $B^I(du)$ an asymmetric quadratic form in $du$.
That is, it must be a linear combination of the gauge-type null
forms
$$Q^I_{JK,jk}(du)=\bigl(\partial_ju^J\partial_ku^K-
\partial_ju^K\partial_ku^J\bigr),
\quad 0\le j<k\le 3, \, 1\le J\le K\le D.$$

It seems likely that one could also allow diagonal terms involving
the relativistic null forms $Q^I_0(du)=(\partial_0u^I)^2-c_I^2
|\nabla_x u^I|^2$, by using a gauge transformation to reduce to
the above types of equations; however, this case will not be
explored here.  One should also be able to allow cubic quasilinear
nonlinearities of the form $R(u,du,d^2u)$ that vanish to second
order in the last two variables.  Doing this, though, would
require handling more powers of $L$, which would complicate the
combinatorics in the continuity argument used to prove global
existence.

In order to solve \eqref{1.5} we must also assume that the data
satisfies the relevant compatibility conditions.  Since these are
well known (see e.g., \cite{KSS}), we shall describe them briefly.
To do so we first let $J_ku =\{\partial^\alpha_xu: \, 0\le
|\alpha|\le k\}$ denote the collection of all spatial derivatives
of $u$ of order up to $k$.  Then if $m$ is fixed and if $u$ is a
formal $H^m$ solution of \eqref{1.5} we can write
$\partial_t^ku(0,\cdot)=\psi_k(J_kf,J_{k-1}g)$, $0\le k\le m$, for
certain compatibility functions $\psi_k$ which depend on the
nonlinear term $Q$ as well as $J_kf$ and $J_{k-1}g$.  Having done
this, the compatibility condition for \eqref{1.5} with $(f,g)\in
H^m\times H^{m-1}$ is just the requirement that the $\psi_k$
vanish on $\partial{\mathcal K}$ when $0\le k\le m-1$.
Additionally, we shall say that $(f,g)\in C^\infty$ satisfy the
compatibility conditions to infinite order if this condition holds
for all $m$.

 We can now state our main result:

\begin{theorem}\label{theorem1.1}
Let ${\mathcal K}$ be a fixed compact obstacle with smooth
boundary that satisfies \eqref{energydecay}.  Assume also that
$Q(du,d^2u)$ and $\square$ are as above. Suppose that $(f,g)\in
C^\infty({\mathbb R}^3\backslash \mathcal{K})$ satisfies the
compatibility conditions to infinite order. Then there is a
constant $\varepsilon_0>0$, and an integer $N > 0$ so that for all
$\varepsilon \leq \varepsilon_0$, if
\begin{equation}\label{1.14}
\sum_{|\alpha|\le N}\|\langle x\rangle^{|\alpha|}
\partial_x^\alpha
 f\|_2 + \sum_{|\alpha|\le N-1}\|
\langle x\rangle^{1+|\alpha|}\partial_x^\alpha g\|_2 \le
\varepsilon,\end{equation}
then \eqref{1.5} has a unique solution
$u\in C^\infty([0,\infty)\times {\mathbb R}^3\backslash
\mathcal{K})$.
\end{theorem}

This paper is organized as follows.  In the next section we shall
collect the $L^2$ estimates that will be needed for the proof of
this existence theorem.  In \S 3 we shall prove the necessary
pointwise decay estimates that will be needed. The results in
these two sections involve variants of those in \cite{KSS3}.  \S 4
will include weighted estimates that are related to the null
condition, which are obstacle variants of ones for the Minkowski
space setting (cf. Hidano \cite{Hidano}, Sideris and Tu
\cite{Si3}, Sogge \cite{So2}, and Yokoyama \cite{Y}). Finally, in
\S 5, we shall use all of these estimates to prove the global
existence theorem.

We are very grateful to S. Zelditch for pointing out the work of
Ikawa and many other suggestions.  It is also a pleasure to thank
N. Burq and S. Klainerman for helpful conversations that
simplified the exposition.  The second author is also grateful for
his collaboration with M. Keel and H. Smith that preceded this
paper.  Both authors would like to thank M. Nakamura for many
helpful comments that greatly helped the exposition.

\newsection{$L^2$ estimates}

We shall be concerned with solutions $u\in
C^\infty(\mathbb{R}_+\times \mathbb{R}^3 \backslash \mathcal{K})$
of the Dirichlet-wave equation
\begin{equation}\label{2.1}
\begin{cases} \square_\gamma u = F
\\
u|_{\partial \kappa}=0
\\
u|_{t=0}=f, \quad \partial_tu|_{t=0}=g,
\end{cases}
\end{equation}
where
$$(\square_\gamma u)^I = (\partial_t^2-c^2_I\Delta)u^I +
\sum_{J=1}^D\sum_{j,k=0}^3 \gamma^{IJ,jk}(t,x)\partial_j\partial_k
u^J, \, \, 1\le I\le D.$$  We shall assume that the
$\gamma^{IJ,jk}$ satisfying the symmetry conditions
\begin{equation}\label{2.2}
\gamma^{IJ,jk}=\gamma^{JI,jk}=\gamma^{IJ,kj}, \end{equation}
 as
well as the size condition
\begin{equation}\label{2.3}
\sum_{I,J=1}^D \sum_{j,k=0}^3\|\gamma^{IJ,jk}(t,x)\|_\infty \le
\delta/(1+t)\end{equation}
 for $\delta>0$ sufficiently small
(depending only on the wave speeds).  The energy estimate will
involve bounds for the gradient of the perturbation terms
$$
\|\gamma'(t,\cd)\|_\infty =\sum_{I,J=1}^D
\sum_{j,k,l=0}^3\|\partial_l\gamma^{IJ,jk}(t,x)\|_\infty,$$ and it
will of course involve the energy form associated with
$\square_\gamma$, $e_0(u)=\sum_{I=1}^D e^I_0(u),$ where
\begin{multline}\label{2.4}
e^I_0(u)= (\partial_0u^I)^2+\sum_{k=1}^3 c^2_I (\partial_ku^I)^2
\\
+2\sum_{J=1}^D\sum_{k=0}^3\gamma^{IJ,0k}\partial_0u^I\partial_ku^J-\sum_{J=1}^D
\sum_{j,k=0}^3\gamma^{IJ,jk}\partial_ju^I\partial_ku^J.
\end{multline}

The most basic estimate will involve
$$E_M(t)=E_M(u)(t)=\int \sum_{j=0}^Me_0(\partial^j_t u)(t,x)\,  dx.$$

\begin{lemma}\label{lemma2.1}  Fix $M=0,1,2,\dots$ and assume
that the perturbation terms $\gamma^{IJ,jk}$ are as above. Suppose
also that if $u\in C^\infty$ solves \eqref{2.1}, and if for every
$t$, $u(t,x)=0$ for large $x$.  Then there is an absolute constant
$C$ so that
\begin{equation}\label{2.5}
\partial_t E^{1/2}_M(t) \le
C\sum_{j=0}^M\|\square_\gamma \partial_t^ju(t,\cd)\|_2 +C\|
\gamma'(t,\cd)\|_\infty E^{1/2}_M(t).
\end{equation}
\end{lemma}

Although the result is standard, we shall present its proof since
it serves as a model for the more difficult variations that are to
follow. We first notice that it suffices to prove the result for
$M=0$ in view of our assumption that the $\partial_t^ju$ satisfy
the Dirichlet boundary conditions for $1\le j\le M$.

To proceed, we need to define the other components of the
energy-momentum vector. For $I=1,2,\dots,D$, and $k=1,2,3$, we let
\begin{equation}\label{2.6}
e^I_k=e^I_k(u)=-2\,c^2_I\,\partial_0u^I\partial_ku^I+2\sum_{J=1}^D\sum_{j=0}^3
\gamma^{IJ,jk}\partial_0u^I\partial_ju^J.
\end{equation}
Then if $e_0$ is the component defined before in \eqref{2.4}, we
have
\begin{align}\label{2.7}
\partial_0e_0^I =&2\,\partial_0u^I\partial_0^2u^I + 2\sum_{k=1}^3
c^2_I\partial_ku^I\partial_0\partial_ku^I +
2\,\partial_0u^I\sum_{J=1}^D\sum_{k=0}^3\gamma^{IJ,0k}\partial_0\partial_ku^J
\\
&+2\sum_{J=1}^D\sum_{k=0}^3
\gamma^{IJ,0k}\partial^2_0u^I\partial_ku^J \notag
\\
&-\sum_{J=1}^D\sum_{j,k=0}^3 \gamma^{IJ,jk}\bigl[
\partial_0\partial_j u^I\partial_ku^J + \partial_j
u^I\partial_0\partial_ku^J \bigr] +R_0^I, \notag
\end{align}
where
$$
R^I_0=2\sum_{J=1}^D \sum_{k=0}^3
(\partial_0\gamma^{IJ,0k})\partial_0u^I\partial_ku^J
-\sum_{J=1}^D\sum_{j,k=0}^3 (\partial_0\gamma^{IJ,jk})\partial_j
u^I\partial_k u^J.
$$
Also,
\begin{align}\label{2.8}
\sum_{k=1}^3\partial_ke^I_k =&-2\,\partial_0u^Ic^2_I\Delta u^I
-2\sum_{k=1}^3c_I^2\partial_ku^I\partial_0\partial_ku^I
\\
& +2\,\partial_0u^I\sum_{J=1}^D\sum_{j=0}^3\sum_{k=1}^3
\gamma^{IJ,jk}\partial_j\partial_ku^J \notag
\\
&+2\sum_{J=1}^D\sum_{j=0}^3\sum_{k=1}^3
\gamma^{IJ,jk}\partial_0\partial_k u^I\partial_ju^J +\sum_{k=1}^3
R^I_k, \notag
\end{align}
where
$$R^I_k=2\sum_{J=1}^D\sum_{j=0}^3
(\partial_k\gamma^{IJ,jk})\partial_0u^I\partial_ju^J.$$

Note that by the symmetry conditions \eqref{2.2} if we sum the
second to last term and the third to last terms in \eqref{2.7}
over $I$, we get
$$-2\sum_{I,J=1}^D \sum_{j=0}^3\sum_{k=1}^3
\gamma^{IJ,jk}\partial_0\partial_ku^I\partial_ju^J,$$ which is
$-1$ times the sum over $I$ of the second to last term of
$\eqref{2.8}$.  From this, we conclude that if we set
$$e_j=e_j(u)=\sum_{I=1}^D e^I_j, \quad j=0,1,2,3,$$
and
$$R=R(u',u')=\sum_{I=1}^D \sum_{k=0}^3 R^I_k,$$
then
$$\partial_t e_0+\sum_{k=1}^3\partial_ke_k  = 2\langle \partial_tu,
\square_\gamma u\rangle
+ R(u',u'),$$ with $\langle \, \cdot \, ,\, \cdot \, \rangle$
denoting the standard inner product in $\R^D$.

If we integrate this identity over $\R^3\backslash\mathcal{K}$ and
apply the divergence theorem, we obtain
\begin{multline}\label{2.9}
\partial_t\int_{\R^3\backslash \mathcal{K}}e_0(t,x) \, dx -
\int_{\partial\mathcal{K}} \sum_{j=1}^3 e_j n_j \, d\sigma
\\
=2\int_{\R^3\backslash \mathcal{K}} \langle
\partial_tu, \square_\gamma u\rangle \, dx
+\int_{\R^3\backslash \mathcal{K}}R(u',u') \, dx\,.
\end{multline}
Here, $n=(n_1,n_2,n_3)$ is the outward normal to $\mathcal{K}$,
and $d\sigma$ is surface measure on $\partial\mathcal{K}$.

Since we are assuming that $u$ solves \eqref{2.1}, and hence
$\partial_t u$ vanishes on $\partial\mathcal{K}$, the integrand in
the last term in the left side of \eqref{2.9} vanishes
identically. Therefore, we have
$$
\partial_t\int_{\R^3\backslash \mathcal{K}} e_0(t,x)\, dx =
2\int_{\R^3\backslash \mathcal{K}} \langle
\partial_tu, \square_\gamma u\rangle\,  dx
+\int_{\R^3\backslash \mathcal{K}}R(u',u')\, dx.
$$
Note that if $\delta$ in \eqref{2.3} is small, then
\begin{equation}
\label{2.10} \bigl(2\max_I\{c^2_I,
c^{-2}_I\}\bigr)^{-1}|u'(t,x)|^2\le e_0(t,x)\le 2\max_I\{c^2_I,
c^{-2}_I\}|u'(t,x)|^2.
\end{equation}
This yields
\begin{multline*}
\partial_t\Bigl(\,\int_{\R^3\backslash
\mathcal{K}} e_0(t,x)\, dx\Bigr)^{1/2}
\\
\le C \|\square_\gamma u(t,\cd)\|_{L^2(\R^3\backslash
\mathcal{K})} +C\sum_{I,J=1}^D\sum_{i,j,k=0}^3\|\partial_i
\gamma^{IJ,jk}(t,\cd)\|_{\infty} \Bigl(\,\int_{\R^3\backslash
\mathcal{K}} e_0(t,x)\,dx\Bigr)^{1/2},
\end{multline*}
as desired. \qed

We require a minor modification of this energy estimate that
involves a slight variant of the scaling vector field
$L=r\partial_r + t\partial_t$.  

Before stating the next result, let us introduce some notation. If
$P=P(t,x,D_t,D_x)$ is a differential operator, we shall let
$$[P,\gamma^{kl} \partial_k\partial_l]u=\sum_{1\le I,J\le D}\sum_{0\le k,l\le 3}
|[P,\gamma^{IJ,kl}\partial_k\partial_l]u^J|.$$
 We can now state the simple variant of Lemma
\ref{lemma2.1} that we require.

\begin{lemma}\label{lemma2.2}
Fix a bump function $\eta\in C^\infty(\mathbb{R}^3)$ satisfying
$\eta(x)=0$, for $x\in \mathcal{K}$  and $\eta(x)=1$, $|x|>1$.
 Let $\tilde L=\eta(x) r\partial_r+t\partial_t$, and set
$$X_{\nu,j}=\int e_0(\tilde L^\nu \partial^j_tu)(t,x)\,
dx.$$ Then if $u\in
C^\infty(\mathbb{R}_+\times\mathbb{R}^3\backslash \mathcal{K})$
solves \eqref{2.1} and vanishes for large $x$ for every $t$
\begin{align}\partial_t X_{\nu,j}&\le CX^{1/2}_{\nu,j}\|\tilde L^\nu
\partial_t^j\square_\gamma
u(t,\cd)\|_2 + C\|\gamma'(t,\cd)\|_\infty X_{\nu,j} \label{2.11}
\\
&+CX^{1/2}_{\nu,j}\|\, [\tilde L^\nu
\partial^j_t, \gamma^{kl} \partial_k\partial_l]u(t,\cd)\|_2 +
CX^{1/2}_{\nu,j}\sum_{\mu\le \nu-1}\|L^\mu
\partial^j_t \square u(t,\cd)\|_2 \notag
\\
&+CX^{1/2}_{\nu,j}\sum_{\substack{\mu+|\alpha|\le j+\nu \\
\mu\le \nu-1}} \|L^\mu \partial^\alpha u'(t,\cd)\|_{L^2(\{x\in
\mathbb{R}^3\backslash \mathcal{K}: \, |x|<1\})}. \notag
\end{align}
\end{lemma}

\noindent {\bf Proof:}  Note that like $u$, $\tilde
L^\nu\partial^j_tu(t,x)$ vanishes when $x\in \partial
\mathcal{K}$.   Therefore by the special case where $M=0$ in Lemma
2.1 we have
\begin{equation}\label{2.12}\parital_t
X_{\nu,j}\le CX_{\nu,j}^{1/2}\|\square_\gamma \tilde L^\nu
\partial_t^j u(t,\cd)\|_2 + C\|\gamma'(t,\cd)\|_\infty X_{\nu,j}.\end{equation}
To proceed we need to estimate the first
term in the right by noting that
\begin{align*}
|\square_\gamma \tilde L^\nu \partial_t^j u|&\le |\tilde L^\nu
\partial_t^j\square_\gamma u|+|[\tilde L^\nu
\partial^j_t, \gamma^{kl} \partial_k\partial_l]u|+|[\tilde
L^\nu , \square ]\partial_t^ju|
\\
&\le |\tilde L^\nu
\partial_t^j\square_\gamma u|+|[\tilde L^\nu
\partial^j_t, \gamma^{kl} \partial_k\partial_l]u| +
|[L^\nu,\square]\partial_t^j u|
\\ &\qquad\qquad\qquad\qquad
+ |[\tilde L^\nu-L^\nu, \square ]\partial^j_tu|
\\
&\le |\tilde L^\nu
\partial_t^j\square_\gamma u|+|[\tilde L^\nu
\partial^j_t, \gamma^{kl} \partial_k\partial_l]u|
+2\sum_{\mu\le \nu-1}|L^\mu\partial^j_t\square u|
\\ &\qquad\qquad\qquad\qquad
+C\chi_{|x|<1}(x)\sum_{\substack{\mu+|\alpha|\le j+\nu \\ \mu\le
\nu-1}}|L^\mu \partial^\alpha u'(t,x)|.
\end{align*}
In the last step we used the fact that $[\square, L]=2\square$,
and $\nabla_x \eta(x)=0$, $|x|>1$.  If we combine the last
inequality and \eqref{2.12} we get \eqref{2.11}. 
 \qed


The last lemma involved estimates for powers of $L$ and
$\partial_t$. Let us now prove a simple result which shows how
these lead to estimates for powers of $L$ and
$\partial=\partial_{t,x}$.

\begin{lemma}\label{lemma2.3}  Fix $N_0$ and $\nu$ and suppose
that $u\in C^\infty(\mathbb{R}_+\times \mathbb{R}\backslash
\mathcal{K})$ solves \eqref{2.1} and vanishes for large $x$ for
each $t$.  Then
\begin{multline}\label{2.13}
\sum_{|\alpha|\le N_0}\|L^\nu \partial^\alpha u'(t,\cd)\|_2\le
C\sum_{\substack{j+\mu \le \nu + N_0 \\ \mu \le \nu}}\|L^\mu
\partial^j_tu'(t,\cd)\|_2
\\
+C\sum_{\substack{|\alpha|+\mu \le N_0+\nu-1 \\ \mu \le
\nu}}\|L^\mu \partial^\alpha \square u(t,\cd)\|_2.\end{multline}
\end{lemma}

\noindent{\bf Proof:}  We shall prove this inequality by induction
on $\nu$.  Since, by elliptic regularity estimates, the inequality
holds when $\nu=0$, let us therefore assume that it is valid when
$\nu$ is replaced by $\nu-1$ and use this to prove it for a given
$\nu=1,2,3,\dots$.

Since $\mathcal{K}\subset \{|x|<1\}$ it is straightforward to see
that
$$\sum_{|\alpha|\le N_0}\|L^\nu \partial^\alpha
u'(t,\cd)\|_{L^2(|x|>1)}$$ is dominated by the right side of
\eqref{2.13}.  Therefore, it suffices to show that we can prove
the analog of \eqref{2.13} where the norm is taken over $|x|<2$.

For the latter, we shall use the fact that
$$\sum_{|\alpha|\le N_0}\|L^\nu \partial^\alpha
u'(t,\cd)\|_{L^2(|x|<2)}\le C\sum_{\substack{|\alpha|+\mu \le
N_0+\nu \\ \mu \le \nu}}t^\mu \|\partial^\mu_t\partial^\alpha
u'(t,\cd)\|_{L^2(|x|<2)}.$$ By elliptic regularity,
\begin{multline*}
\sum_{\substack{|\alpha|+\mu \le N_0+\nu \\ \mu \le
\nu}}\|\partial^\alpha \partial^\mu_t u'(t,\cd)\|_{L^2(|x|<2)} \le
C\sum_{\substack{j+\mu \le N_0+\nu \\ \mu \le
\nu}}\|\partial_t^{j+\mu}u'(t,\cd)\|_{L^2(|x|<4)}
\\
+C\sum_{\substack{|\alpha|+\mu \le N_0+\nu -1 \\ \mu\le
\nu}}\|\partial^\alpha \partial^\mu_t \square
u(t,\cd)\|_{L^2(|x|<4)}.
\end{multline*}
Therefore,
\begin{align*}
&\sum_{|\alpha|\le N_0}\|L^\nu \partial^\alpha
u'(t,\cd)\|_{L^2(|x|<2)}
\\
&\le C\sum_{\substack{j+\mu \le N_0+\nu \\ \mu\le \nu}}\|t^\mu
\partial_t^{\mu+j}u'(t,\cd)\|_{L^2(|x|<4)}+C\sum_{\substack{|\alpha|+\mu\le
N_0+\nu-1 \\ \mu\le \nu}}\|t^\mu \partial_t^\mu \partial^\alpha
\square u(t,\cd)\|_{L^2(|x|<4)}
\\
&\le C\sum_{j\le N_0}\|L^\nu
\partial_t^j u'(t,\cd)\|_2 + C\sum_{\substack{|\alpha|+\mu \le
N_0+\nu \\ \mu\le \nu-1}}\|L^\mu \partial^\alpha u'(t,\cd)\|_2
\\
&\qquad\qquad\qquad\qquad\qquad\qquad
+C\sum_{\substack{|\alpha|+\mu\le N_0+\nu-1 \\
\mu\le \nu}}\|L^\mu
\partial^\alpha \square u(t,\cd)\|_2.
\end{align*}
As a result, we get \eqref{2.13} by the inductive step and the
fact that, we can control the norms over the set where $|x|>1$.
\qed

Using \eqref{2.13} we can prove the following estimate.

\begin{proposition}\label{proposition2.4}  Suppose that the constant $\delta$
in \eqref{2.3} is small.  Suppose further that
\begin{equation}\label{2.14} \|\gamma'(t,\cd)\|_\infty \le
\delta/(1+t),
\end{equation}
and
\begin{multline}\label{2.15}
\sum_{\substack{j+\mu \le N_0+\nu_0 \\ \mu\le
\nu_0}}\Bigl(\|\tilde L^\mu \partial^j_t \square_\gamma
u(t,\cd)\|_2 +\|\, [\tilde L^\mu \partial^j_t,
\gamma^{kl}\partial_k\partial_l]u(t,\cd)\|_2\Bigr)
\\
\le \frac{\delta}{1+t}\sum_{\substack{j+\mu \le N_0+\nu_0 \\
\mu\le \nu_0}}\|\tilde L^\mu \partial^j_t u'(t,\cd)\|_2 +
H_{\nu_0,N_0}(t),
\end{multline}
where $N_0$ and $\nu_0$ are fixed.  Then
\begin{align}\label{2.16}&\sum_{\substack{|\alpha|+\mu\le
N_0+\nu_0\\ \mu\le \nu_0}}\|L^\mu\partial^\alpha u'(t,\cd)\|_2
\\
&\le C\sum_{\substack{|\alpha|+\mu\le N_0+\nu_0-1 \\ \mu\le
\nu_0}}\|L^\mu \partial^\alpha \square
u(t,\cd)\|_2+C(1+t)^{A\delta} \sum_{\substack{\mu+j\le N_0+\nu_0
\\ \mu\le \nu_0}} X^{1/2}_{\mu,j}(0)
\notag
\\
& +C(1+t)^{A\delta}\Bigl(\int_0^t \sum_{\substack{|\alpha|+\mu \le
N_0+\nu_0-1 \\ \mu\le \nu_0-1}}\|L^\mu \partial^\alpha \square
u(s,\cd)\|_2\, ds+\int_0^t H_{\nu_0,N_0}(s)ds\Bigr) \notag
\\
&+C(1+t)^{A\delta}\int_0^t \sum_{\substack{|\alpha|+\mu \le
N_0+\nu_0 \\ \mu \le \nu_0-1}}\|L^\mu \partial^\alpha
u'(s,\cd)\|_{L^2(|x|<1)}\, ds, \notag
\end{align}
where the constants $C$ and $A$ depend only on the constants in
\eqref{2.11}.
\end{proposition}

In practice $H_{\nu_0,N_0}(t)$ will involve weighted $L^2_x$ norms
of $|L^\mu\partial^\alpha u'|^2$ with $\mu+|\alpha|$ much smaller
than $N_0+\nu_0$, and so the integral involving $H_{\nu_0,N_0}$
can be dealt with using an inductive argument and weighted
$L^2_tL^2_x$ estimates that will be presented at the end of this
section.

\noindent{\bf Proof:}  We first note that by \eqref{2.3} and the
definition \eqref{2.4} of the energy form
\begin{equation}\label{2.17}
\sum_{\substack{j+\mu \le N_0+\nu_0\\ \mu\le \nu_0}}\|\tilde L^\mu
\partial^j_t u'(t,\cd)\|_2 \le 2\sum_{\substack{j+\mu\le N_0+\nu_0
\\ \mu\le \nu_0}}X^{1/2}_{\mu,j}(t),
\end{equation}
if $\delta$ is sufficiently small.
  Therefore, by \eqref{2.11} and
\eqref{2.14}-\eqref{2.15} we have
\begin{align*}\partial_t \sum_{\substack{j+\mu \le N_0+\nu_0 \\
\mu\le \nu_0}}X^{1/2}_{\mu,j}(t)&\le
\frac{A\delta}{1+t}\sum_{\substack{j+\mu \le N_0+\nu_0 \\
\mu\le \nu_0}}X^{1/2}_{\mu,j}(t) + AH_{\nu_0,N_0}(t)
\\
&+A\sum_{\substack{\mu+j\le N_0+\nu_0-1\\ \mu\le \nu_0-1}}\|L^\mu
\partial_t^j \square u(t,\cd)\|_2
\\
&+A\sum_{\substack{|\alpha|+\mu \le N_0+\nu_0 \\ \mu\le
\nu_0-1}}\|L^\mu \partial^\alpha u'(t,\cd)\|_{L^2(|x|<1)},
\end{align*}
where $A$ depends on the constants in \eqref{2.11}.  By Gronwall's
inequality
$$\sum_{\substack{j+\mu  \le N_0+\nu_0 \\ \mu\le \nu_0}}
X^{1/2}_{\mu,j}$$ is dominated by the right side of \eqref{2.16}.
 By applying \eqref{2.13} and \eqref{2.17}, we conclude that
\eqref{2.16} must be valid.  \qed

nonnegative

In proving our existence results for \eqref{1.5} the key step will
be to obtain a priori $L^2$-estimates involving $L^\mu Z^\alpha
u'$.  The next result indicates how these can be obtained from
ones involving $L^\mu\partial^\alpha u'$.

\begin{proposition}\label{proposition2.5}  Fix $N_0$ and $\nu_0$,
and set
\begin{equation}\label{2.18}
Y_{N_0,\nu_0}(t)=\sum_{\substack{|\alpha|+\mu\le
N_0+\nu_0 \\ \mu\le \nu_0}} \int e_0(L^\mu Z^\alpha u)(t,x)\,
dx.\end{equation} Suppose that the constant $\delta$ in
\eqref{2.3} is small and that \eqref{2.14} holds.  Then
\begin{multline}\label{2.19}
\partial_t Y_{N_0,\nu_0} \le C Y^{1/2}_{N_0,\nu_0} \sum_{\substack{
|\alpha|+\mu\le N_0+\nu_0\\\mu\le\nu_0}} \|\Box_\gamma L^\mu Z^\alpha
u(t,\cd)\|_2 + C\|\gamma'(t,\cd)\|_\infty Y_{N_0,\nu_0} \\
+C \sum_{\substack{|\alpha|+\mu\le N_0+\nu_0+1\\ \mu\le \nu_0}}
\|L^\mu \partial^\alpha u'(s,\cd)\|^2_{L^2(|x|<1)}.
\end{multline}
\end{proposition}

\noindent{\bf Proof:}  If we argue as in the proof of Lemma
\ref{lemma2.1} we find that
\begin{multline}\label{2.20}\partial_t Y_{N_0,\nu_0}\le
CY^{1/2}_{N_0,\nu_0}\sum_{\substack{|\alpha|+\mu \le N_0+\nu_0 \\
\mu\le \nu_0}}\|\square_\gamma L^\mu Z^\alpha u(t,\cd)\|_2 +
C\|\gamma'(t,\cd)\|_\infty Y_{N_0,\nu_0}
\\
+C\int_{\partial \mathcal{K}}\sum_{a=1}^3 |e_an_a|\, d\sigma,
\end{multline}
where $n=(n_1,n_2,n_3)$ is the outward normal at a given point in
$\partial \mathcal{K}$, and
$$e_a=\sum_{\substack{|\alpha|+\mu \le N_0+\nu_0 \\ \mu\le
\nu_0}}e_a(L^\mu Z^\alpha u)(t,x), \quad a=1,2,3,$$ are the
components of the energy-momentum tensor defined in \eqref{2.6}.
Since $\mathcal{K}\subset \{|x|<1\}$ and since
$$\sum_{\substack{|\alpha|+\mu \le N_0+\nu_0 \\ \mu\le
\nu_0}}|L^\mu Z^\alpha u(t,x)|\le C\sum_{\substack{|\alpha|+\mu
\le N_0+\nu_0 \\ \mu\le \nu_0}}|L^\mu\partial^\alpha u(t,x)|,
\quad x\in \partial \mathcal{K},$$ we have
$$\int_{\partial \mathcal{K}}\sum_{a=1}^3 |e_an_a|\, d\sigma\le
C\int_{\{x\in \mathbb{R}^3\backslash \mathcal{K}: \, |x|<1\}}
\sum_{\substack{|\alpha|+\mu \le N_0+\nu_0+1 \\ \mu\le
\nu_0}}|L^\mu \partial^\alpha u'(t,x)|^2 \, dx.$$  \qed

As in \cite{KSS2} and \cite{KSS3} we shall control the local $L^2$
norms, such as the last term in \eqref{2.19} by using weighted
$L^2_tL^2_x$ estimates.  They will also be used in obtaining decay
estimate for solutions of the nonlinear equation. To avoid
cumbersome notation, for the rest of the section we shall abuse
notation a bit by letting $\square = \partial^2_t -\Delta$ denote
the unit speed D'Alembertian.  The passing from the ensuing
estimates involving this case to ones involving \eqref{1.6} is
straightforward.  Also, in what follows, we shall let
$$S_T=\{[0,T]\times \mathbb{R}^3\backslash \mathcal{K}\}$$
denote the time strip of height $T$ in $\mathbb{R}_+\times
\mathbb{R}^3\backslash \mathcal{K}$.

\begin{proposition}\label{proposition2.6}  Fix $N_0$ and $\nu_0$.
Suppose that $\mathcal{K}$ satisfies the local exponential energy
decay bounds \eqref{energydecay}.  Suppose also that $u\in
C^\infty$ solves \eqref{2.1} and satisfies $u(t,x)=0$, $t<0$. Then
there is a constant $C=C_{N_0,\nu_0,\mathcal{K}}$ so that if $u$
vanishes for large $x$ for every fixed $t$
\begin{align}\label{2.21}
\bigl(\log (2+T)\bigr)^{-1/2}\sum_{\substack{|\alpha|+\mu \le
N_0+\nu_0 \\ \mu\le \nu_0}}&\|\langle x \rangle^{-1/2}L^\mu
\partial^\alpha u'\|_{L^2(S_T)}
\\
&\le C\int_0^T\sum_{\substack{|\alpha|+\mu \le N_0+\nu_0+1 \\
\mu\le \nu_0}}\|\square L^\mu \partial^\alpha u(s,\cd)\|_2\, ds
\notag
\\
&+C\sum_{\substack{|\alpha|+\mu \le N_0+\nu_0 \\
\mu\le \nu_0}}\|\square L^\mu \partial^\alpha u\|_{L^2(S_T)}
\notag
\end{align}
Also, if $N_0$ and $\nu_0$ are fixed
\begin{align}\label{2.22}
\bigl(\log (2+T)\bigr)^{-1/2}\sum_{\substack{|\alpha|+\mu \le
N_0+\nu_0 \\ \mu\le \nu_0}}&\|\langle x \rangle^{-1/2}L^\mu
Z^\alpha u'\|_{L^2(S_T)}
\\
&\le C\int_0^T\sum_{\substack{|\alpha|+\mu \le N_0+\nu_0+1 \\
\mu\le \nu_0}}\|\square L^\mu Z^\alpha u(s,\cd)\|_2\, ds \notag
\\
&+C\sum_{\substack{|\alpha|+\mu \le N_0+\nu_0 \\
\mu\le \nu_0}}\|\square L^\mu Z^\alpha u\|_{L^2(S_T)}. \notag
\end{align}
\end{proposition}

To prove these estimates we shall need a couple of lemmas. The
first says that these estimates hold (with no loss of derivatives)
in the boundaryless case.

\begin{lemma}\label{lemma2.7}  Suppose that $v\in
C^\infty(\mathbb{R}\times \mathbb{R}^3)$ vanishes for large $x$
for every $t$. Then there is a uniform constant $C$ so that if $v$
has vanishing Cauchy data
\begin{equation}\label{2.23}
(\ln(2+T))^{-1/2}\|\langle x \rangle^{-1/2} v'\|_{L^2([0,T]\times
\mathbb{R}^3)} \le
C\int_0^T\|\square v(s,\cd)\|_{L^2(\mathbb{R}^3)}\, ds.
\end{equation}
Also, given $\mu$ and $\alpha$,
\begin{multline}\label{2.24}
(\ln(2+T))^{-1/2}\|\langle x\rangle^{-1/2}(L^\mu Z^\alpha
v)'\|_{L^2([0,T]\times \mathbb{R}^3)}
\\
\le
C\int_0^T\|\square L^\mu Z^\alpha v(s,\cd)\|_{L^2(\mathbb{R}^3)}\,
ds.
\end{multline}
\end{lemma}

The first inequality, \eqref{2.23}, was proved in \cite{KSS2}. The
second follows from the first.

As was shown in \cite{KSS2}, \eqref{2.23} follows immediately from
the fact that stronger bounds hold when one restricts the norms in
the left to regions where $|x|$ is bounded.  In particular, just
by using Huygens principle, one can show that if $R$ is fixed then
there is a uniform constant $C=C_R$ so that
\begin{equation}\label{2.25}
\|v'\|_{L^2([0,T]\times \{x\in \mathbb{R}^3: |x|<R\})}
+\|v\|_{L^2_tL^6_x([0,T]\times \{x\in \mathbb{R}^3: |x|<R\})}
\le
C\int_0^T\|\square
v(s,\cd)\|_{L^2(\mathbb{R}^3)}\, ds.
\end{equation}

To prove Proposition \ref{proposition2.6} we shall need the
following local estimates which follow from the local exponential
energy decay \eqref{energydecay}.

\begin{lemma}\label{lemma2.8}  Suppose that  \eqref{energydecay}
holds and that $\square u(t,x)=0$ for $|x|>4$ and $t>0$. Suppose also that
$u(t,x)=0$ for $t\le 0$.  Then if $N_0$ and $\nu_0$ are fixed and if
$c>0$ is as in \eqref{energydecay}
\begin{align}\label{2.26}
\sum_{\substack{|\alpha|+\mu \le N_0+\nu_0 \\ \mu\le
\nu_0}}\|L^\mu \partial^\alpha
u'(t,\cd)&\|_{L^2(\{\mathbb{R}^3\backslash \mathcal{K}: \,
|x|<4\})}
\\
&\le C\sum_{\substack{|\alpha|+\mu\le N_0+\nu_0-1 \\
\mu\le \nu_0}}\|L^\mu \partial^\alpha \square u(t,\cd)\|_2
\notag
\\
&+C\int_0^te^{-(c/2)(t-s)}\sum_{\substack{|\alpha|+\mu \le N_0+\nu_0+1
\\ \mu\le \nu_0}}\|L^\mu \partial^\alpha \square
u(s,\cd)\|_2 ds \notag 
\end{align}
\end{lemma}

The proof is quite simple.  By \eqref{energydecay} we have that
$$\sum_{\substack{j+\mu \le N_0+\nu_0 \\ \mu\le \nu_0}}\|\langle
t\rangle^\mu \partial_t^\mu
\partial_t^ju'(t,\cd)\|_{L^2(\{\mathbb{R}^3\backslash \mathcal{K}:
\, |x|<6\})}$$ is dominated by the last term
in the right side of \eqref{2.26}.  By Lemma \ref{lemma2.3}, this
implies that \eqref{2.26} holds.

\noindent{\bf Proof of Proposition \ref{proposition2.6}:}  We
shall only prove \eqref{2.21} since \eqref{2.22} follows from the
same argument.

The first step in proving \eqref{2.21} will be to show that if we
take the $L^2_tL^2_x$ norm over a region where $|x|$ is bounded
then we have better estimates, i.e.,
\begin{align}\label{2.27}
\sum_{\substack{|\alpha|+\mu\le N_0+\nu_0 \\ \mu\le \nu_0}}\|L^\mu
\partial^\alpha u'\|_{L^2(S_T\cap |x|<2)}
&\le C\int_0^T\sum_{\substack{|\alpha|+\mu\le N_0+\nu_0+1 \\
\mu\le \nu_0}}\|\square L^\mu\partial^\alpha u(s,\cd)\|_2ds
\\
&+C\sum_{\substack{|\alpha|+\mu\le N_0+\nu_0 -1\\ \mu\le
\nu_0}}\|\square L^\mu \partial^\alpha u\|_{L^2(S_T)} \notag
\end{align}

To prove this, let us first assume that $u$ is as in Lemma
\ref{lemma2.8}.  Thus, if we assume that $\square u(t,x)=0$
when $|x|>4$, then by \eqref{2.26} we have for $0<\tau<T$
\begin{align*}
\sum_{\substack{|\alpha|+\mu \le N_0+\nu_0 \\ \mu\le
\nu_0}}&\|L^\mu \partial^\alpha u'(\tau,\cd)\|^2_{L^2(|x|<2)}
\\
&\le C\sum_{\substack{|\alpha|+\mu \le N_0+\nu_0-1 \\ \mu\le
\nu_0}}\|L^\mu \partial^\alpha \square
u(\tau,\cd)\|_2^2
\\
&+C\Bigl(\int_0^\tau e^{-c(\tau-s)}\sum_{\substack{|\alpha|+\mu
\le N_0+\nu_0+1 \\ \mu\le \nu_0}}\|L^\mu\partial^\alpha \square
u(s,\cd)\|_2ds\Bigr)
\\
&\qquad\qquad\qquad\times\Bigl(\int_0^\tau\sum_{\substack{|\alpha|+\mu
\le N_0+\nu_0+1
\\ \mu\le \nu_0}}\|L^\mu\partial^\alpha \square
u(s,\cd)\|_2ds\Bigr).
\end{align*}
After integrating $\tau$ from $0$ to $T$ we obtain \eqref{2.27}
under the support assumptions of Lemma \ref{lemma2.8}.

Note that if we had applied Young's inequality, we would have
gotten
\begin{multline}\label{2.28}
\sum_{\substack{|\alpha|+\mu \le N_0+\nu_0 \\ \mu\le
\nu_0}}\|L^\mu\partial^\alpha u'\|^2_{L^2(S_T\cap |x|<2)}
\le C\sum_{\substack{|\alpha|+\mu \le N_0+\nu_0+1 \\ \mu\le
\nu_0}}\|L^\mu\partial^\alpha \square u\|^2_{L^2(S_T)}
,
\\
\text{if } \square u(t,x)=u(0,x)=\partial_tu(0,x)=0 \, , \, |x|>4.
\end{multline}
This inequality will be useful in the last part of the proof of
\eqref{2.27} where we need to show that the inequality holds when
we assume that 
$\square u(t,x)$ vanishes for $|x|<3$.  To do this, we fix
$\rho\in C^\infty(\mathbb{R}^3)$ satisfying $\rho(x)=1$ for
$|x|<2$ and $\rho(x)=0$ for $|x|\ge 3$.  We then write $u=u_0+u_r$
where $u_0$ solves the boundaryless wave equation $\square
u_0(t,x)=\square u(t,x)$ if $|x|\ge 3$ and $0$ otherwise with
vanishing
initial data 
. It then follows that $\tilde u=\rho u_0+u_r$ solves the
Dirichlet-wave equation $\square \tilde
u=-2\nabla_x\rho\cdot\nabla_xu_0 - (\Delta \rho)u_0$ with zero
initial data.  Therefore, by \eqref{2.28}, we have
\begin{align*}
\sum_{\substack{|\alpha|+\mu \le N_0+\nu_0 \\ \mu\le
\nu_0}}&\|L^\mu\partial^\alpha u'\|^2_{L^2(S_T\cap
|x|<2)}=\sum_{\substack{|\alpha|+\mu \le N_0+\nu_0 \\ \mu\le
\nu_0}}\|L^\mu\partial^\alpha \tilde u'\|^2_{L^2(S_T\cap |x|<2)}
\\
&\le \sum_{\substack{|\alpha|+\mu \le N_0+\nu_0+1 \\ \mu\le
\nu_0}}\|L^\mu\partial^\alpha \square\tilde u\|^2_{L^2(S_T)}
\\
&\le C\sum_{\substack{|\alpha|+\mu \le N_0+\nu_0 +1\\ \mu\le
\nu_0}}\Bigl(\|L^\mu\partial^\alpha u'_0\|^2_{L^2(S_T\cap |x|<4)}
+\|L^\mu\partial^\alpha u_0\|^2_{L^2(S_T\cap |x|<4)}\Bigr).
\end{align*}
One now gets \eqref{2.27} for this by applying \eqref{2.25} since
$\square u_0=\square u$ in $\mathbb{R}^3\backslash \mathcal{K}$
.

To finish the proof of \eqref{2.21} we must show that
\begin{align}\label{2.29}
\bigl(\log &(2+T)\bigr)^{-1/2}\sum_{\substack{|\alpha|+\mu \le
N_0+\nu_0 \\ \mu\le \nu_0}}\|\langle x \rangle^{-1/2}L^\mu
\partial^\alpha u'\|_{L^2(S_T\cap |x|>2)}
\\
&
\le C\int_0^T\sum_{\substack{|\alpha|+\mu \le N_0+\nu_0+1 \\
\mu\le \nu_0}}\|\square L^\mu \partial^\alpha u(s,\cd)\|_2\, ds
%
+C\sum_{\substack{|\alpha|+\mu \le N_0+\nu_0-1 \\
\mu\le \nu_0}}\|\square L^\mu \partial^\alpha u\|_{L^2(S_T)}
\notag
\end{align}

To do this, we fix $\beta\in C^\infty(\mathbb{R}^3)$ satisfying
$\beta(x)=1$, $|x|\ge 2$ and $\beta(x)=0$, $|x|\le 3/2$.  By
assumption the obstacle is contained in the set $|x|<1$.  It
follows that $v=\beta u$ solves the boundaryless wave equation
$\square v=\beta \square u-2\nabla_x\beta\cdot
\nabla_xu-(\Delta\beta)u$ with vanishing initial data.
Also $u(t,x)=v(t,x)$ for $|x|>2$.  We split $v=v_1+v_2$ where
$v_1$ solves $\square v_1=\beta \square u$
and $v_2$ solves $\square v_2=-2\nabla_x\beta\cdot
\nabla_xu-(\Delta \beta)u$
and both have
zero initial data. By \eqref{2.24} if we replace $u$ by $v_1$ in
the left side of \eqref{2.29}, then the resulting quantity is
dominated by the right side of \eqref{2.29}.

Therefore, to finish the proof, we must show that
\begin{multline}\label{2.30}
(\log(2+T))^{-1/2}\sum_{\substack{|\alpha|+\mu\le N_0+\nu_0\\
\mu\le \nu_0}}\|\langle x\rangle^{-1/2}L^\mu\partial^\alpha
v_2'\|_{L^2(S_T\cap |x|>2)}
\\
\le C\int_0^T \sum_{\substack{|\alpha|+\mu\le N_0+\nu_0+1 \\
\mu\le \nu_0}}\|\square L^\mu \partial^\alpha u(s,\cd)\|_2 ds
+C\sum_{\substack{|\alpha|+\mu\le N_0+\nu_0-1\\
\mu\le \nu_0}}\|\square L^\mu\partial^\alpha u\|_{L^2(S_T)}.
\end{multline}
To prove this, we note that
$G=-2\nabla_x\beta\cdot\nabla_xu-(\Delta\beta)u=\square v_2$
vanishes unless $1<|x|<2$.  To use this, fix $\chi\in
C^\infty_0(\mathbb{R})$ satisfying $\chi(s)=0$ for $s\not\in [1/2, 2]$ and
$\sum_j\chi(s-j)=1$.  We then split $G=\sum_jG_j$ where
$G_j(s,x)=\chi(s-j)G(s,x)$, and let $v_{2,j}$ be the solution of
the inhomogeneous wave equation $\square v_{2,j}=G_j$ in Minkowski
space with zero initial data.  Since $v_2$ also has vanishing
Cauchy data, by the sharp Huygens principle the functions
$v_{2,j}$ have finite overlap, so that we have
$|L^\mu\partial^\alpha v_2'|^2\le C\sum_j|L^\mu\partial^\alpha
v_{2,j}'|^2$ for some uniform constant $C$.  Therefore, by
\eqref{2.24}, the square of the left side of \eqref{2.30} is
dominated by
\begin{align*}
\sum_{\substack{|\alpha|+\mu\le N_0+\nu_0 \\ \mu\le
\nu_0}}&\sum_j\Bigl(\int_0^T\|L^\mu \partial^\alpha
G_j(s,\cd)\|_{L^2(\mathbb{R}^3)} \, ds \Bigr)^2
\\
&\le C\sum_{\substack{|\alpha|+\mu\le N_0+\nu_0 \\
\mu\le \nu_0}} \|L^\mu\partial^\alpha G\|^2_{L^2([0,T]\times
\mathbb{R}^3)}
\\
&\le C\sum_{\substack{|\alpha|+\mu\le N_0+\nu_0 \\ \mu\le
\nu_0}}\Bigl(\|L^\mu\partial^\alpha u'\|^2_{L^2([0,T]\times
|x|<2)} + \|L^\mu \partial^\alpha u\|^2_{L^2([0,T]\times
|x|<2)}\Bigr)
\\
&\le \sum_{\substack{|\alpha|+\mu\le N_0+\nu_0 \\ \mu\le
\nu_0}}\|L^\mu \partial^\alpha u'\|^2_{L^2(S_T\cap |x|<2)}.
\end{align*}
Consequently, the bound \eqref{2.30} follows from \eqref{2.27}.

This finishes the proof of \eqref{2.21}.  Since the other part of
the proposition follow from the same argument, this completes the
proof of Proposition \ref{proposition2.6}.  \qed

To be able to handle the last term in the right side of
\eqref{2.16} we shall need the following result which follows from
a similar argument.

\begin{lemma}\label{lemma2.9}
Suppose that \eqref{energydecay} holds, and suppose that $u\in
C^\infty$ solves \eqref{2.1} and satisfies $u(t,x)=0$, $t<0$.
Then if $\nu_0$ and $N_0$ are fixed and if $c>0$ is as in
\eqref{energydecay}
\begin{align}\label{2.31}
&\sum_{\substack{|\alpha|+\mu\le N_0+\nu_0 \\ \mu\le
\nu_0}}\|L^\mu
\partial^\alpha u'(t,\cd)\|_{L^2(|x|<2)}
\\
&\le C\sum_{\substack{|\alpha|+\mu\le N_0+\nu_0+1 \\ \mu\le
\nu_0}}\Biggl[\int_0^t e^{-\frac{c}{2}(t-s)}\|L^\mu
\partial^\alpha \square u(s,\cd)\|_{L^2(|x|<4)} \, ds
+\|L^\mu\partial^\alpha \square u(t,\cd)\|_{L^2(|x|<4)}\Biggr]
\notag
\\
&+C\sum_{\substack{|\alpha|+\mu\le N_0+ \nu_0+1 \\ \mu\le
\nu_0}}\int_0^t e^{-\frac{c}{2}(t-s)}\Bigl(\int_0^s\|L^\mu \partial^\alpha
\square u(\tau,\cd)\|_{L^2(| \, |x|-(s-\tau)\, |<10)} \,
d\tau\Bigr) ds \notag
\\
&+C\sum_{\substack{|\alpha|+\mu\le N_0+ \nu_0+1 \\ \mu\le
\nu_0}}\int_0^t\|L^\mu \partial^\alpha \square u(s,\cd)\|_{L^2(|\,
|x|-(t-s)\, |<10)} \, ds .\notag
\end{align}
Additionally, if $t>2$,
\begin{multline}\label{2.32}
\sum_{\substack{|\alpha|+\mu\le N_0+\nu_0 \\ \mu\le
\nu_0}}\int_0^t\|L^\mu \partial^\alpha u'(s,\cd)\|_{L^2(|x|<2)}\,
ds
\\
\le C\sum_{\substack{|\alpha|+\mu\le N_0+\nu_0+1 \\ \mu\le
\nu_0}}\int_0^t \Biggl(\int_0^s \|L^\mu \partial^\alpha \square
u(\tau,\cd)\|_{L^2(|\, |x|-(s-\tau)\, |<10)}\, d\tau\Biggr)\, ds.
\end{multline}
\end{lemma}

\noindent {\bf Proof:}  Since the first inequality obviously
implies the second, we shall only prove \eqref{2.31}.

If $\square u(s,x)$ vanishes when $|x|>4$, the result follows from
\eqref{2.26}.  In this case a stronger inequality holds where the
last term in the right is not present.

To finish we need to show that the inequality is valid when
$\square u(s,x)$ vanishes for $|x|<3$.  In this case, as in the
proof of Proposition \ref{proposition2.6} we write $u=u_0+u_r$
where $u_0$ solves $\square u_0=\square u$ with vanishing Cauchy
data.  Then if as above $\rho\in C^\infty(\mathbb{R}^3)$ equals
$1$ for $|x|<2$ and $0$ for $|x|>3$, then $\tilde u=\rho u_0+u_r$
has vanishing Cauchy data and solves $\square \tilde
u=-2\nabla_x\rho \cdot \nabla_x u_0 - (\Delta \rho)u_0$. Thus,
since $\square \tilde u=0$ for $|x|>3$, by the above case
\begin{align*}
&\sum_{\substack{|\alpha|+\mu \le N_0+\nu_0 \\ \mu\le
\nu_0}}\|L^\mu \partial^\alpha u'(t,\cd)\|_{L^2(|x|<2)}
\\
&= \sum_{\substack{|\alpha|+\mu \le N_0+\nu_0 \\ \mu\le
\nu_0}}\|L^\mu \partial^\alpha \tilde u'(t,\cd)\|_{L^2(|x|<2)}
\\
&\le C\sum_{\substack{|\alpha|+\mu \le N_0+\nu_0+1 \\
\mu\le \nu_0}}\Bigl[\int_0^t e^{-\frac{c}{2}(t-s)}\|L^\mu \partial^\alpha
\square \tilde u(s,\cd)\|_2 ds+ \|L^\mu\partial^\alpha \square
\tilde u(t,\cd)\|_2\Bigr]
\\
&\le C\sum_{\substack{|\alpha|+\mu \le N_0+\nu_0+1 \\
\mu\le \nu_0}}\Bigl[\int_0^t e^{-\frac{c}{2}(t-s)}\bigl(\|L^\mu
\partial^\alpha u'_0(s,\cd)\|_{L^2(|x|<4)}+\|L^\mu \partial^\alpha
u_0(s,\cd)\|_{L^2(|x|<4)}\bigr)ds
\\
&\qquad \qquad\qquad\qquad\qquad+\|L^\mu\partial^\alpha
u'_0(t,\cd)\|_{L^2(|x|<4)}+\|L^\mu
\partial^\alpha u_0(t,\cd)\|_{L^2(|x|<4)}\Bigr] .
\end{align*}
Since $\square u=\square u_0$ one can use the sharp Huygens
principle to see that the last term is dominated by the last term
in the right side of \eqref{2.31}, which finishes the proof.  \qed

\newsection{Pointwise Estimates}

We will estimate solutions of the scalar inhomogeneous wave
equation
\begin{equation}\label{3.1}
\begin{cases}
(\partial_t^2-\Delta) w(t,x)=F(t,x), \quad (t,x)\in
\mathbb{R}_+\times \mathbb{R}^3\backslash \mathcal{K}
\\
w(t,x)=0, \quad x\in \partial\mathcal{K}
\\
w(t,x)=0, \quad t\le 0.
\end{cases}\end{equation}

If we assume, as before, we assume that $\mathcal{K}\subset \{x\in
\mathbb{R}^3: \, |x|<1\}$ then we have the following

\begin{theorem}\label{theorem3.1}  Suppose that the local
energy decay bounds \eqref{energydecay} hold for
$\mathcal{K}$.  Suppose also that $|\alpha|=M$.  Then
\begin{align}\label{3.2}
(1+t+|x|) |L^\nu Z^\alpha w(t,x)| &\le C\int_0^t \int_{\R^3
\backslash \mathcal{K}}\sum_{\substack{|\beta|+\mu\le M+\nu+7\\
\mu\le \nu+1}} |L^\mu Z^\beta F(s,y)|\frac{dyds}{|y|}
\\
&+C\int_0^t \sum_{\substack{|\beta|+\mu\le M+\nu +4 \\ \mu\le \nu
+1}}\|L^\mu
\partial^\beta F(s,\cd)\|_{L^2(\{x\in \R^3 \backslash
\mathcal{K}: \, |x|<2\})} \, ds. \notag
\end{align}
\end{theorem}

The special case of this estimate where $\nu=0$ was handled in
\cite{KSS3} in the non-trapping case.  Since it is technically
harder to handle pointwise bounds involving powers of $L$, we
shall give the proof of \eqref{3.2} for the sake of completeness.
Handling the case where there is a loss of regularity in the
energy decay as in \eqref{energydecay} does not present any added
difficulty.  The fact that \eqref{energydecay} involves a loss of
one derivative accounts why when $\nu=0$ the right side of
\eqref{3.2} involves on extra derivative versus the results in
\cite{KSS3}.

The proof will resemble that of Proposition \ref{proposition2.6}.
We shall prove the estimate when $x$ is near the obstacle
primarily by using the local energy decay estimates
\eqref{energydecay}, while away from the obstacle we shall mainly
use the fact that related bounds hold in Minkowski space.

The Minkowski space estimates we shall use say that if $w_0$ is a
solution of the inhomogeneous wave equation
\begin{equation}\label{3.3}
\begin{cases} (\partial_t^2-\Delta)w_0(t,x)=G(t,x), \quad (t,x)\in
\R_+\times \R^3
\\
w_0(0,x)=\partial_tw_0(0,x)=0,
\end{cases}
\end{equation}
then
\begin{equation}\label{3.4}
(1+t+|x|)|L^\nu Z^\alpha w_0(t,x)|\le C\int_0^t
\int_{\R^3}\sum_{\substack{|\beta|+\mu \le |\alpha|+\nu+3\\
\mu\le \nu+1}}|L^{\mu}Z^\beta G(s,y)|\, \frac{dyds}{|y|}.
\end{equation}
This follows from inequalities (2.3) and (2.9) in \cite{KSS3} and
the fact that $[\partial^2_t-\Delta, Z]=0$, and
$[\partial_t^2-\Delta, L]=2(\partial^2_t-\Delta)$.  The estimate
where the weight in the left is $(1+t)$ was the main pointwise
estimate in \cite{KSS3}, while the contribution of the weight
$|x|$ in the left just follows from the fact that
\begin{align}\label{3.5}
|x|\, |w_0(t,x)|&\le C\int_0^t\int_{|\, |x|-(t-s)\, |}^{|x|+(t-s)}
\sup_{|\theta|=1}|G(s,r\theta)|\, rdr ds
\\
&\le C\int_0^t\int_{\{y\in \R^3: \, |y|\in [|\, |x|-(t-s)\, |,
|x|+(t-s)]\}}\sum_{|a|\le 2}|\Omega^a G(s,y)| \, \frac{dyds}{|y|}.
\notag
\end{align}

Recall that we are assuming that $\mathcal{K}\subset \{x\in \R^3:
\, |x|<1\}$.  With this in mind, the first step is to see that
\eqref{3.4} and \eqref{3.5} yield
\begin{align}\label{3.6}
(1+t+|x|)|L^\nu Z^\alpha w(t,x)|&\le \int_0^t \int_{\R^3\backslash
\mathcal{K}}\sum_{\substack{|\beta|+\mu\le |\alpha|+\nu+3\\ \mu\le\nu +
1}}|L^{\mu}Z^{\beta}F(s,y)|\, \frac{dyds}{|y|}
\\
&+C\sup_{|y|\le 2, 0\le s\le t}(1+s)\bigl( |L^\nu Z^\alpha
w'(s,y)|+|L^\nu Z^\alpha w(s,y)|\bigr).
\notag
\end{align}
The proof is exactly like that of Lemma 4.2 in \cite{KSS3}.  One
fixes $\rho\in C^\infty(\R)$ satisfying $\rho(r)=1$, $r\ge 2$,
$\rho(r)=0$, $r\le 1$, and then applies \eqref{3.4}-\eqref{3.5} to
$w_0(t,x)=\rho(|x|)L^\nu Z^\alpha w(t,x)$, which solves the
inhomogeneous wave equation
\begin{multline*}
(\partial^2_t-\Delta)w_0(t,x)=\rho(|x|)(\partial_t^2-\Delta)L^\nu
Z^\alpha w(t,x)
\\
-2\rho'(|x|)\frac{x}{|x|}\cdot \nabla_x L^\nu Z^\alpha
w(t,x)-(\Delta \rho(|x|))L^\nu Z^\alpha w(t,x),
\end{multline*}
with zero initial data.  When one applies \eqref{3.4}, the first
term in the right side of this equation results in the first term
in the right side of of \eqref{3.6}, while if one applies the
first inequality in \eqref{3.5} one sees that the last two terms
of the equation result in the last two terms of \eqref{3.6}.

It remains to prove pointwise bounds in the region where $|x|<2$.
Additionally since the coefficients of $Z$ are bounded, it
suffices to show that if $|\gamma|\le |\alpha|+1=M+1$, then
\begin{align}\label{3.7}
t\sup_{|x|<2}|L^\nu \partial^\gamma w(t,x)|&\le
C\int_0^t\int_{\R^3\backslash
\mathcal{K}}\sum_{\substack{|\beta|+\mu \le M+\nu+7 \\ \mu\le
\nu+1}} |L^\mu Z^\beta F(s,y)|\, \frac{dyds}{|y|}
\\
&+C\int_0^t\sum_{\substack{|\beta|+\mu\le M+\nu+4 \\ \mu\le
\nu+1}}\|L^\mu
\partial^\beta F(s,\cdot)\|_{L^2(\{x\in \R\backslash
\mathcal{K}: \, |x|<4\})} \, ds.
\notag
\end{align}

Using cutoffs for the forcing terms, we can split things into
proving \eqref{3.7} for the following two cases
\begin{itemize}
\item {\bf Case 1:} $F(s,y)=0$ if $|y|>4$
\item {\bf Case 2:} $F(s,y)=0$ if $|y|<3$.
\end{itemize}

For either case, we shall use the following immediate consequence
of the Fundamental Theorem of Calculus:
$$|tL^\nu \partial^\gamma w(t,x)|\le \sum_{j=0,1}
\int_0^t|(s\partial_s)^jL^\nu\partial^\gamma w(s,x)|\, ds.$$ If we
apply the Sobolev lemma, using the fact that $|\gamma|\le M+1$,
and that Dirichlet conditions allow us to control $w$ locally by
$w'$, then we get
\begin{align*}
t\sup_{|x|<2}|L^\nu \partial^\gamma w(t,x)|&\le C\int_0^t
\sum_{|\beta|\le M+2, \mu\le 1}\|(s\partial_s)^\mu
L^\nu\partial^\beta w'(s,\cdot)\|_{L^2(\R^3\backslash \mathcal{K}:
\, |x|<4)} \, ds
\\
&\le C\int_0^t\sum_{\substack{|\beta|+\mu\le M+\nu+3 \\ \mu\le
\nu+1}} \|L^\mu\partial^\beta w'(s,\cd)\|_{L^2(\R^3\backslash
\mathcal{K}: \, |x|<4)} \, ds.
\end{align*}
If we are in Case 1, we apply \eqref{2.26} to get the variant of
\eqref{3.7} involving only the second term in the right.

In Case 2, we need to write $w=w_0+w_r$ where $w_0$ solves the
boundaryless wave equation $(\partial_t^2-\Delta)w_0=F$ with zero
initial data.  Fix $\eta\in C^\infty_0(\R^3)$ satisfying
$\eta(x)=1$, $|x|<2$ and $\eta(x)=0$, $|x|\ge 3$.  It then follows
that if we set $\tilde w=\eta w_0 +w_r$, then since $\eta F=0$,
$\tilde w$ solves the Dirichlet-wave equation
$$(\partial_t^2-\Delta)\tilde w=G=-2\nabla_x\eta \cdot
\nabla_xw_0-(\Delta \eta)w_0$$ with zero initial data.  The
forcing term vanishes unless $2\le |x|\le 4$.  Hence, by Case 1
\begin{align*}
t\sup_{|x|<2}|L^\nu\partial^\gamma
w(t,x)|&=t\sup_{|x|<2}|L^\nu\partial^\gamma\tilde w(t,x)|
\\
&\le C\int_0^t \sum_{\substack{|\beta|+\mu\le M+\nu+4 \\
\mu\le \nu+1}} \|L^\mu\partial^\beta
G(s,\cd)\|_{L^2(\R^3\backslash \mathcal{K})}\, ds
\\
&\le C\int_0^t \sum_{\substack{|\beta|+\mu\le M+\nu+5 \\
\mu\le \nu+1}}\|L^\mu\partial^\beta w_0(s,\cd)\|_{L^\infty(2\le
|x|\le 4)}\, ds.
\end{align*}

To finish the argument, we apply \eqref{3.5} to obtain
$$\|L^\mu\partial^\beta w_0(s,\cd)\|_{L^\infty(2\le |x|\le 4)}
\le C\sum_{|a|\le 2}\int_0^s\int_{|s-\tau-|y|\, |\le 4}
|L^\mu\partial^\beta \Omega^a F(\tau,y)|\, \frac{dyd\tau}{|y|}.$$
Note that the sets $\Lambda_s=\{(\tau,y): \, 0\le \tau \le s,
|s-\tau-|y|\, |\le 4\}$ satisfy $\Lambda_s\cap
\Lambda_{s'}=\emptyset$ if $|s-s'|\ge 10$.  Therefore, if in the
preceding inequality we sum over $|\beta|+\mu\le M+\nu+5$, $\mu\le
\nu+1$, and then integrate over $s\in [0,t]$ we conclude that
\eqref{3.7} must hold for Case 2, which finishes the proof. \qed

%
\eqref{energydecay} holds.

\newsection{Estimates related to the null condition}

Here we shall prove simple bounds for the null forms.  They must
involve the weight $<c_Jt-r>$ due to the fact that we are not
using the generators of Lorentz rotations.  The estimates will
involve the admissible homogeneous vector fields that we are using
$\{\Gamma\}=\{Z, L\}$.  Also, as before, $\partial$ denotes the
space-time gradient $\nabla_{t,x}$.

\begin{lemma}\label{lemma4.1}  Suppose that the quasilinear
null condition \eqref{1.9} holds.  Then
\begin{multline}\label{4.1}
\Bigl|\sum_{0\le j,k,l\le 3}B^{IJ,jk}_{Jl}\partial_l
u\partial_j\partial_k v\Bigr|
\\
\le C<r>^{-1}\bigl( |\Gamma u|\,
|\partial^2 v| +|\partial u|\, |\partial\Gamma v|\bigr) +
C\frac{<c_Jt-r>}{<t+r>} |\partial u|\, |\partial^2 v|.
\end{multline}
Also, if the asymmetric semilinear null condition \eqref{1.10}
holds
\begin{equation}\label{4.2}
\Bigl|\sum_{0\le j,k\le 3}A^{I,jk}_{JK}\partial_j u\partial_k
v\Bigr| \le C<r>^{-1}\bigl(|\Gamma u|\, |\partial v| + |\partial
u|\, |\Gamma v|\bigr)
.
\end{equation}
\end{lemma}

\noindent{\bf Proof:}  The first estimate is well known.  See,
e.g., \cite{Si3}, \cite{So2}.  It also follows from the proof of
\eqref{4.2}.

Proving \eqref{4.2} is straightforward.  Since we are assuming
\eqref{1.10} the quadratic form involved must be skew symmetric.
If we write $\nabla_x = \frac{x}r \partial_r + \frac{x}{r^2}
\wedge \Omega$, then since $|(\frac{x}{r^2}\wedge \Omega)u|\le
C\langle r \rangle^{-1}|\Gamma u|$, we conclude that the left side
of \eqref{4.2} must be dominated by
$$\langle r\rangle^{-1}\bigl(|\Omega u|\, |\partial v|+|\partial
u|\, |\Omega v|\bigr) + |\partial_t u\partial_rv-\partial_ru
\partial_tv|.$$ If we write $\partial_r=r^{-1}L
+\frac{t}{r}\partial_t$ then we can estimate the last term
$$|\partial_tu\partial_rv-\partial_ru\partial_tv|\le
\tfrac1r\bigl(|Lu|\, |\partial_tv|+|\partial_tu|\, |Lv|\bigr).$$
Combining these two steps yields \eqref{4.2}.  \qed

\medskip

We also need the following result.

\begin{lemma}\label{lemma4.2}  If $h\in
C^\infty_0$ has Dirichlet boundary conditions then if $R<t/2$ and
$t\ge 1$
\begin{align}\label{4.3}
\|\partial &h'(t,\cd)\|_{L^2(R/2<|x|<R)}
\\
&\le Ct^{-1}\bigl(\sum_{|\alpha|\le 1}\|\Gamma^\alpha
h'(t,\cd)\|_{L^2(R/4<|x|<2R)}
 + t\|(\partial_t^2-\Delta)
h(t,\cd)\|_{L^2(R/4<|x|<2R)}\bigr)
\notag
\\
&
+C\|\langle x\rangle^{-1}h'(t,\cd)\|_{L^2(R/4<|x|<2R)} +
C\|\langle x\rangle^{-2}h(t,\cd)\|_{L^2(R/4<|x|<2R)}. \notag
\end{align}
Also,
\begin{multline}\label{4.4}
\|<t-r>\partial h'(t,\cd)\|_{L^2(|x|>t/4)}
\\
\le C\sum_{|\alpha|\le 1}\|\Gamma^\alpha h'(t,\cd)\|_2 +
C\|<t+r>(\partial_t^2-\Delta) h(t,\cd)\|_2,
\end{multline}
and if $\delta>0$ is fixed then
\begin{multline}\label{4.5}
\|h'(t,\cd)\|_{L^6(|x|\notin [(1-\delta)t,(1+\delta)t], \,
|x|>\delta t)}
\\
\le Ct^{-1}\Bigl(\sum_{|\alpha|\le 1}\|\Gamma^\alpha h'(t,\cd)\|_2
+ C\|<t+r>(\partial_t^2-\Delta) h(t,\cd)\|_2\Bigr).
\end{multline}
\end{lemma}

\noindent{\bf Proof:} To prove \eqref{4.3} we need to use the fact
(see \cite{KS}, Lemma 2.3) that
\begin{equation}\label{4.6}
\langle t-r\rangle \bigl(|\partial \partial_t h(t,x)|+|\Delta
h(t,x)|\bigr) \le C\sum_{|\alpha|\le 1}|\partial\Gamma^\alpha
h(t,x)|+C\langle t+r\rangle |(\partial^2_t-\Delta)h(t,x)|.
\end{equation}
Also, elliptic regularity gives
\begin{multline*}\|\nabla_x
h'(t,\cd)\|_{L^2(|x|\in [R/2,R])}\le C\|\Delta
h(t,\cd)\|_{L^2(|x|\in [R/4,2R])}
\\
+CR^{-1}\|h'(t,\cd)\|_{L^2( |x|\in
[R/4,2R])}+CR^{-2}\|h(t,\cd)\|_{L^2( |x|\in [R/4,2R])}.
\end{multline*}
If we combine these two inequalities then we get \eqref{4.3}

 To
prove \eqref{4.4} we need to use another estimate from \cite{KS},
namely, if $g\in C^\infty_0(\R_+\times \R^3)$,
$$\|\langle t-r\rangle \nabla^2_x g(t,\cd)\|_{L^2(\R^3)}\le
C\sum_{|\alpha|\le 1}\|\Gamma^\alpha
g'(t,\cd)\|_{L^2(\R^3)}+C\|\langle t+r\rangle
(\partial^2_t-\Delta)g\|_{L^2(\R^3)}.$$ If we fix $\eta\in
C^\infty(\R^3)$ satisfying $\eta(x)=1$, $|x|>1/4$ and $\eta(x)=0$,
$|x|<1/8$ and let $g(t,x)=\eta(x/\langle t\rangle)h(t,x)$ then we
conclude that the analog of \eqref{4.4} must hold where $\nabla
h'$ is replaced by $\nabla_x h'$.  Since \eqref{4.6} yields the
same bounds for $\partial_t h'$, we get \eqref{4.4}.

Inequality \eqref{4.5} follows from the fact that its left side is
dominated by
$$\|\nabla_x h'(t,\cd)\|_{L^2(|x|\notin
[(1-\delta/2)t,(1+\delta/2)t], \, |x|>\delta t/2)} +
t^{-1}\|h'(t,\cd)\|_2.$$ Since the proof of \eqref{4.4} implies
that the first term is dominated by the right side of \eqref{4.5}
if $\delta>0$ is fixed, we are done. \qed

The following result will be useful for dealing with waves
interacting at different speeds.

\begin{corr}\label{corra.3} Fix $c_1, c_2>0$ satisfying $c_1\ne c_2$.  Then
if
$u,v\in C^\infty_0(\R_+\times \R^3\backslash
\mathcal{K})$ vanish on $\R_+\times \partial\mathcal{K}$
\begin{align}\label{4.7}
&\int_{\mathbb{R}^3\backslash\mathcal{K}}|\partial^2 u(t,x)|\,
|v'(t,x)| <x>^{-1} dx
\\
& \le Ct^{-1}\Bigl(\sum_{|\alpha|\le 1}\|\Gamma^\alpha
u'(t,\cd)\|_2 + \|<t+r>(\partial_t^2-c_1^2\Delta)
u(t,\cd)\|_2\Bigr) \|<x>^{-1} v'(t,\cd)\|_2 \notag
\\
& +C \sum_{R=2^k\le t/2} \bigl(\|\langle
x\rangle^{-1}u'(t,\cd)\|_{L^2( R/2<|x|<R)}
+
\|\langle x \rangle^{-2}u(t,\cd)\|_{L^2(R/2<|x|<R)} \bigr) \notag
\\
&\qquad\qquad \times \|<x>^{-1}v'(t,\cd)\|_{L^2( R/2<|x|<R)}
\notag
\\
&+Ct^{-4/3}\bigl(\sum_{|\alpha|\le 1}\|\Gamma^\alpha u'(t,\cd)\|_2
+ \|<t+r>(\partial_t^2-c_1^2\Delta) u(t,\cd)\|_2\bigr) \notag
\\
&\qquad\qquad\times
 \bigl(\sum_{|\alpha|\le
1}\|\Gamma^\alpha v'(t,\cd)\|_2+\|<t+r>(\partial_t^2-c_2^2\Delta)
v(t,\cd)\|_2\bigr) \notag
\end{align}
\end{corr}

\noindent{\bf Proof:}  Let $\delta<|c_1-c_2|$.  Then if we use
 Schwarz's inequality, \eqref{4.3} and \eqref{4.4} we see that we can bound
 $$\int_{|x|\notin ((1-\delta)c_1t,(1+\delta)c_1t)}|\partial^2 u(t,x)|\,
 |v'(t,x)| <x>^{-1} dx$$
 by the first two terms in the right side of \eqref{4.7}.

For a given $j=0,1,2,\dots$ we can use H\"older's inequality, to
find that
\begin{multline*}\int_{<c_1t-r>\in [2^j,2^{j+1})} |\partial^2 u(t,x)| \,
|v'(t,x)|
<x>^{-1} dx
\\
\le C t^{-1/3}2^{j/3}\|\partial^2 u(t,\cd)\|_{L^2(<c_1t-r>\in
(2^j, 2^{j+1}))}\|v'(t,\cd)\|_{L^6(<c_1t-r>\in (2^j, 2^{j+1}))},
\end{multline*}
assuming that $r$ is bounded below by a fixed multiple of $t$ when
$<c_1t-r>\in [2^j,2^{j+1})$. Since $\delta<|c_1-c_2|$, if $\{x: \,
<c_1t-r>\in [2^j,2^{j+1})\} \cap \{x: r\in
((1-\delta)c_1t,(1+\delta)c_1t)\} \ne \emptyset$, we can apply
\eqref{4.4} and \eqref{4.5} to see that the right side is bounded
by $2^{-2j/3}$ times the second term in the right side of
\eqref{4.7}. After summing over $j$, this implies that when we
restrict the integration in the left side of \eqref{4.7} to the
the set where $r\in ((1-\delta)c_1t, (1+\delta)c_1t)$, the
resulting expression is dominated by the second term in the right
of \eqref{4.7}.  This completes the proof. \qed

 To handle same-speed interactions, we shall need the
following similar result.

\begin{corr}\label{corra.4} Let
$u,v\in C^\infty_0(\R_+\times \R^3\backslash \mathcal{K})$
vanish on $\R_+\times \partial\mathcal{K}$.
 Then,
\begin{align}\label{4.8}
\int_{\mathbb{R}^3\backslash\mathcal{K}}&\frac{\langle
t-r\rangle}{\langle t+r\rangle}|\partial^2 u(t,x)|\, |v'(t,x)|\,
\langle x\rangle^{-1} \, dx
\\
 &\le
Ct^{-1}\bigl(\sum_{|\alpha|\le1}\|\Gamma^\alpha
u'(t,\cd)\|_2+\|\langle t+r\rangle\square
u(t,\cd)\|_2\bigr)\|\langle x\rangle^{-1}v'(t,\cd)\|_2 \notag
\\
&+C\sum_{R=2^k< t/2}\bigl(\|\langle x\rangle^{-1}u'(t,\cd)\|_{L^2(
R/4<|x|<2R)}+
\|\langle x\rangle^{-2}u(t,\cd)\|_{L^2(R/4<|x|<2R)}\bigr) \notag
\\
 &\qquad\qquad \qquad\qquad \times \|\langle
x\rangle^{-1}v'(t,\cd)\|_{L^2( R/4<|x|<2R)}.\notag
\end{align}
\end{corr}

\noindent{\bf Proof of Corollary \ref{corra.4}:}  

To prove \eqref{4.8} we just use   Schwarz's inequality and
\eqref{4.3} and \eqref{4.4} to see that its left side is dominated
by
\begin{align*}
t^{-1}&\|\langle
t-r\rangle\partial^2u(t,\cd)\|_{L^2(|x|>t/4)}\|\langle
x\rangle^{-1}v'(t,\cd)\|_{L^2(|x|>t/4)}
\\
&+ \sum_{R=2^k<t/2}  t^{-1}\|\langle
t-r\rangle\partial^2u(t,\cd)\|_{L^2( R/2<|x|<R)}\|\langle
x\rangle^{-1}v'(t,\cd)\|_{L^2( R/2<|x|<R)}
\\
&\le C t^{-1}\bigl(\sum_{|\alpha|\le1}\|\Gamma^\alpha
u'(t,\cd)\|_2+\|\langle t+r\rangle\square
u(t,\cd)\|_2\bigr)\|\langle x\rangle^{-1}v'(t,\cd)\|_2
\\
&+C\sum_{R=2^k<t/2}\bigl(\|\langle x\rangle^{-1}u'(t,\cd)\|_{L^2(
R/4<|x|<2R)}+
\|\langle x\rangle^{-2}u(t,\cd)\|_{L^2(R/4<|x|<2R)}\bigr)
\\
&\qquad\qquad\qquad\qquad\qquad \times\|\langle
x\rangle^{-1}v'(t,\cd)\|_{L^2(R/4<|x|<2R)},
\end{align*}
 which completes the
proof. \qed

We also need the following consequence of the Sobolev lemma (see
\cite{knull}).

\begin{lemma}\label{lemma4.5}  Suppose that $h\in C^\infty(\mathbb{R}^3)$.  Then
for $R\ge 1$
$$\|h\|_{L^\infty(R/2<|x|<R)}\le CR^{-1}\sum_{|\alpha|+|\beta|\le
2}\|\Omega^\alpha \partial_x^\beta h\|_{L^2(R/4<|x|<2R)}.$$ Also,
$$\|h\|_{L^\infty(R<|x|<R+1)}\le CR^{-1}\sum_{|\alpha|+|\beta|\le
2}\|\Omega^\alpha \partial_x^\beta h\|_{L^2(R-1<|x|<R+2)}.$$
\end{lemma}

\newsection{Continuity argument}

In this section we shall prove our main result, Theorem
\ref{theorem1.1}.  We shall take $N=101$ in its smallness
hypothesis \eqref{1.14}, but this certainly is not optimal.

We start out with a number of straightforward reductions that will
allow us to use the estimates from \S 2--4.

First, let us assume
 that the wave speeds $c_I$ all are distinct
since straightforward modifications of the argument give the more
general case where the various components are allowed to have the
same speed.

To prove our global existence theorem we shall need a standard
local existence theorem:

\begin{theorem}\label{theorem5.1}  Suppose that $f$ and $g$ are as
in Theorem \ref{theorem1.1} with $N\ge6$ in \eqref{1.14}.  Then
there is a $T>0$ so that the initial value problem \eqref{1.5}
with this initial data has a $C^2$ solution satisfying
$$
u\in L^\infty([0,T]; H^{N}(\mathbb{R}^3\backslash
\mathcal{K}))\cap C^{0,1}([0,T]; H^{N-1}(\mathbb{R}^3\backslash
\mathcal{K})).$$ The supremum of such $T$ is equal to the supremum
of all $T$ such that the initial value problem has a $C^2$
solution with $\partial^\alpha u$ bounded for $|\alpha|\le 2$.
Also, one can take $T\ge 2$ if $\|f\|_{H^{N}}+\|g\|_{H^{N-1}}$ is
sufficiently small.
\end{theorem}

This essentially follows from the local existence results Theorem
9.4 and Lemma 9.6 in \cite{KSS}.  The latter were only stated for
diagonal single-speed systems; however, since the proof relied
only on energy estimates, it extends to the multi-speed
non-diagonal case if the symmetry assumptions \eqref{1.8} are
satisfied.

Next, as in \cite{KSS3}, in order to avoid dealing with
compatibility conditions for the Cauchy data, it is convenient to
reduce the Cauchy problem \eqref{1.5} to an equivalent equation
with a nonlinear driving force but vanishing Cauchy data.  We then
can set up a continuity argument for the new equation using the
estimates from \S 2--4 to prove Theorem \ref{theorem1.1}.

Recall that our smallness condition on the data is
\begin{equation}\label{5.1} \sum_{|\alpha|\le 101}\|\langle
x\rangle^\alpha \partial_x^\alpha f\|_{L^2(\mathbb{R}^3\backslash
\mathcal{K})}+\sum_{|\alpha|\le 100}\|\langle
x\rangle^{1+|\alpha|}\partial_x g\|_{L^2(\mathbb{R}^3\backslash
\mathcal{K})}\le \varepsilon.
\end{equation}
To make the reduction to an equation with zero initial data, we
first note that if the data satisfies \eqref{5.1} with
$\varepsilon>0$ small, then we can find a solution $u$ to the
system \eqref{1.5} on a set of the form $0<ct<|x|$ where
$c=5\max_I c_I$, and that this solution satisfies
\begin{equation}\label{5.2}
\sup_{0<t<\infty}\sum_{|\alpha|\le 101}\|\langle x\rangle^\alpha
\partial^\alpha u(t,\cd)\|_{L^2(\mathbb{R}^3\backslash
\mathcal{K}: \, |x|>ct)}\le C_0\varepsilon,
\end{equation}
where $C_0$ is an absolute constant.

To prove this we shall repeat an argument from \cite{KSS3}.  We
note that by scaling in the $t$-variable we may assume that
$\max_Ic_I=1/2$.  The above local existence theorem yields a
solution $u$ to \eqref{1.5} on the set $0<t<2$ satisfying the
bounds \eqref{5.2}.  To see that this solution extends to the
larger set $0<ct<|x|$, we let $R\ge 4$ and consider data
$(f_R,g_R)$ supported in the set $R/4<|x|<4R$ which agrees with
the data $(f,g)$ on the set $R/2<|x|<2R$.  Let $u_R(t,x)$ satisfy
the boundaryless equation
$$\square u_R= Q(du_R,R^{-1}d^2u_R)$$
with Cauchy data $(f_R(R\cdot),Rg_R(R\cdot))$.  The solution $u_R$
then exists for $0<t<1$ by standard results (see \cite{H}) and
satisfies
\begin{multline*}
\sup_{0<t<1}\|u_R(t,\cd)\|_{H^{101}(\mathbb{R}^3)}\le
C\bigl(\|f_R(R\cdot)\|_{H^{101}(\mathbb{R}^3)}+R\|g_R(R\cdot)\|_{H^{100}(\mathbb
{R}^3)}
\\
\le CR^{-3/2}\Bigl(\sum_{|\alpha|\le 101}\|(R\partial_x)^\alpha
f_R\|_{L^2(\mathbb{R}^3)}+R\sum_{|\alpha|\le
100}\|(R\partial_x)^\alpha g_R\|_{L^2(\mathbb{R}^3)}\Bigr).
\end{multline*}
The smallness condition on $|u_R'|$ implies that the wave speeds
for the quasilinear equation are bounded above by 1.  A domain of
dependence argument shows that the solutions
$u_R(R^{-1}t,R^{-1}x)$ restricted to $|\, |x|-R\, |<\tfrac{R}2-t$
agree on their overlaps, and also with the local solution,
yielding a solution to \eqref{1.5} on the set
$\{\mathbb{R}^3\backslash \mathcal{K}: \, 2t<|x|\}$.  A partition
of unity argument now yields \eqref{5.2}.

We use the local solution $u$ to set up the continuity argument.
Fix a cutoff function $\chi\in C^\infty(\mathbb{R})$ satisfying
$\chi(s)=1$ if $s\le \frac{1}{2c}$ and $\chi(s)=0$ if
$s>\frac{1}{c}$, and set
$$u_0(t,x)=\eta(t,x)u(t,x), \quad \quad
\eta(t,x)=\chi(|x|^{-1}t),$$ assuming as we may that
$0\in\mathcal{K}$. Note that since $|x|$ is bounded below on the
complement of $\mathcal{K}$, the function $\eta(t,x)$ is smooth
and homogeneous of degree $0$ in $(t,x)$. Also
$$\square u_0=\eta Q(du,d^2u)+[\square, \eta] u.$$
Thus, $u$ solves $\square u=Q(du,d^2u)$ for $0<t<T$ if and only if
$w=u-u_0$ solves
\begin{equation}\label{5.3}
\begin{cases}\square w=(1-\eta)Q(
du,d^2u)-[\square ,\eta]u
\\
w|_{\partial\mathcal{K}}=0
\\
w(t,x)=0, \quad t\le 0
\end{cases}
\end{equation}
for $0<t<T$.

The key step in proving that \eqref{5.3} admits a global solution
is to prove uniform dispersive estimates for $w$ on intervals of
existence.  To do this, let us first note that since $u_0=\eta u$
by \eqref{5.2} and Lemma \ref{lemma4.5} there is an absolute
constant $C_1$ so that
\begin{multline}\label{5.4}
(1+t+|x|)\sum_{\mu+|\alpha|\le 99}|L^\mu Z^\alpha
u_0(t,x)|
\\ + \sum_{\mu+|\alpha|+|\beta|\le 101}\|\langle t+r\rangle^{|\beta|} L^\mu
Z^\alpha \partial^\beta u_0(t,\cd)\|_2
\le C_1\varepsilon.\end{multline} Furthermore, if we let $v$ be
the solution of the linear equation
\begin{equation}\label{5.5}
\begin{cases}\square v=-[\square,\eta]u
\\
v|_{\partial \mathcal{K}}=0
\\
v(t,x)=0, \quad t\le 0,
\end{cases}
\end{equation}
then \eqref{5.2} and Theorem 3.1 implies that there is an absolute
constant $C_2$ so that
\begin{equation}\label{5.6}
(1+t+|x|)\sum_{\mu+|\alpha|\le 90}|L^\mu Z^\alpha v(t,x)|\le
C_2\varepsilon.
\end{equation}
Indeed, by \eqref{3.2} the left side of \eqref{5.6} is dominated
by
\begin{multline*}
\int_0^t\int_{|x|>cs} \sum_{\mu+|\alpha|\le 97}|L^\mu Z^\alpha
([\square, \eta]u)(s,x)| \, \frac{dx ds}{|x|} \\+ \int_0^t
\sum_{\mu+|\alpha|\le
94}\|L^\mu\partial^\beta([\square,\eta]u)(s,\cdot)\|_{L^2(\mathbb
R\backslash\mathcal K\, :\, |x|<2)}\:ds
\end{multline*}
 which by the Schwarz
inequality is bounded by
$$\sum_{\mu+|\alpha|\le 97}\sum_{j=0}^\infty
\sup_{0<cs<2^{j+1}}\|\langle x\rangle^{3/2} L^\mu Z^\alpha
[\square, \eta] u(s,\cd)\|_{L^2(\mathbb{R}^3\backslash
\mathcal{K}: \, 2^j<|x|<2^{j+1})}.$$ Since this is bounded by
$$\sup_{0<t<\infty}\sum_{\mu+|\alpha|\le 97}\|\langle x\rangle^2
L^\mu Z^\alpha [\square, \eta]u(t,\cd)\|_2,$$ one gets \eqref{5.6}
by \eqref{5.2} and the homogeneity of $\eta$.

Using this we can set up the continuity argument.  If
$\varepsilon>0$ is as above we shall assume that we have a $C^2$
solution of our equation \eqref{1.5} for $0\le t\le T$ such that
for $t\in [0,T]$ and small $\varepsilon>0$ we have the pointwise
dispersive estimates
\begin{align}\label{5.7}
(1+t+r)\sum_{|\alpha| \le 40 }\Bigl(| Z^\alpha w(t,x)|+|Z^\alpha
w'(t,x)|\Bigr)&\le A_0\varepsilon
\\
(1+t+r)\sum_{\substack{|\alpha|+\nu \le 55 \\ \nu\le 2}}|L^\nu
Z^\alpha w(t,x)|&\le B_1\varepsilon(1+t)^{1/5}\log(2+t),
\label{5.8}
\end{align}
as well as the $L^2_x$ and weighted $L^2_tL^2_x$ estimates
\begin{align}\label{5.9}
\sum_{|\alpha|\le 100}\|\partial^\alpha w'(t,\cd)\|_2
&\le
B_2\varepsilon(1+t)^{1/20}
\\
\sum_{\substack{|\alpha|+\nu\le 70\\ \nu\le 3}}\|L^\nu Z^\alpha
w'(t,\cd)\|_2 &\le B_3\varepsilon (1+t)^{1/10} \label{5.10}
\\
\sum_{\substack{|\alpha|+\nu\le 68 \\ \nu\le 3}}\|\langle
x\rangle^{-1/2}L^\nu Z^\alpha w'\|_{L^2(S_t)}&\le B_4\varepsilon
(1+t)^{1/10}(\log(2+t))^{1/2}. \label{5.11}
\end{align}
Here, as before the $L^2_x$-norms are taken over
$\mathbb{R}^3\backslash \mathcal{K}$, and the weighted
$L^2_tL^2_x$-norms are taken over $S_t = [0,t]\times
\mathbb{R}^3\backslash \mathcal{K}$.

In our main estimate, \eqref{5.7}, $A_0=4C_2$, where $C_2$ is the
constant occurring for the bounds \eqref{5.6} for $v$.  Clearly if
$\varepsilon$ is small then all of these estimates are valid if
$T=2$, by Theorem \ref{theorem5.1}. Keeping this in mind, we shall
then prove that for $\varepsilon>0$ smaller than some number
depending on the constants $B_1$-$B_4$ that

{\bf{i)}} \eqref{5.7} is valid with $A_0$ replaced by $A_0/2$;

{\bf{ii)}} \eqref{5.8}--\eqref{5.11} are a consequence of
\eqref{5.7}
for suitable constants $B_i$.

By the local existence theorem it will follow that a solution
exists for all $t>0$ if $\varepsilon$ is small enough.

Let us first deal with i).   Since we already know that $v$
satisfies \eqref{5.6} to achieve i), by Theorem \ref{theorem3.1}
it suffices to show that
\begin{equation}\label{5.12}
I+II \le C\varepsilon^2,
\end{equation}
where
\begin{align} \label{5.13}
I&=\int_0^t \int_{\mathbb{R}^3\backslash
\mathcal{K}}\sum_{\substack{|\alpha|+\mu \le 48 \\ \mu\le
1}}|L^\mu Z^\alpha Q(du,d^2u)(s,y)|\frac{dsdy}{|y|}
\\
II&=\int_0^t \sum_{\substack{|\alpha|+\mu \le 45 \\ \mu\le
1}}\|L^\mu \partial^\alpha Q(du,d^2u)(s,\cd)\|_{L^2( |x|<2)} \,
ds, \label{5.14}
\end{align}
since this implies the same sort of bounds where $Q$ is replaced
by $(1-\eta)Q$ in \eqref{5.13} and \eqref{5.14}.

Let us first deal with $I$.  This term was the only one that had
to be dealt with in the boundaryless case, and the argument for it
is similar to the corresponding one in \cite{So2}.

To handle $I$ we shall have to employ a different argument for the
quadratic terms satisfying the null condition and the quasilinear
ones that do not.  Therefore, let us write
\begin{equation}\label{5.15}
Q=\square u=N(u',u'')+\sum_{J\ne K}\sum_{j,k,l=0}^3
B_{K,l}^{IJ,jk}\partial_lu^K\partial_j\partial_ku^J, \quad
u=u_0+w,
\end{equation}
where the ``null term" $N(u',u'')$ satisfies the bounds in Lemma
\ref{lemma4.1}, while the second term in the right of
\eqref{5.15} involves interactions between waves of different
speeds.

Let us first handle the contribution of $N(u',u'')$ to $I$.  By
Lemma \ref{lemma4.1}
\begin{align}\label{5.16}
\sum_{\substack{|\alpha|+\mu \le 48 \\ \mu\le 1}}|L^\mu Z^\alpha
N(u',u'')| &\le \frac{C}{|y|}\sum_{\substack{|\alpha|+\mu \le 50 \\
\mu\le 2}}|L^\mu Z^\alpha u|\sum_{\substack{|\alpha|+\mu \le 50 \\
\mu\le 2}}|L^\mu Z^\alpha u'|
\\
&+C\sum_J \frac{\langle c_Jt-r\rangle}{\langle
t+r\rangle}\sum_{\substack{|\alpha|+\mu \le 48 \\ \mu\le 1}}|L^\mu
Z^\alpha \partial u|\sum_{\substack{|\alpha|+\mu \le 48 \\ \mu\le
1}}|L^\mu Z^\alpha \partial^2 u|. \notag
\end{align}

To handle the contribution of the first term in the right side of
\eqref{5.16} to $I$, we apply
\eqref{5.8} to get that
$$\sum_{\substack{|\alpha|+\mu \le 50 \\ \mu\le 2}}|L^\mu Z^\alpha
u(y,s)|\le C\varepsilon (|y|+s)^{-4/5}\log (2+s),$$ which means
that the first term in the right side of \eqref{5.16} has a
contribution to $I$ which is dominated by
\begin{multline*}
\varepsilon\int_0^t \int_{\mathbb{R}^3\backslash
\mathcal{K}}\sum_{\substack{|\alpha|+\mu \le 50 \\ \mu\le
2}}|L^\mu Z^\alpha u'(s,y)|\frac{\log(2+s) \, dy ds}{|y|^2
(|y|+s)^{4/5}}
\\
\le C_\delta \varepsilon \int_0^t \sum_{\substack{|\alpha|+\mu \le
50 \\ \mu\le 2}}\|\langle y\rangle^{-1/2}L^\mu Z^\alpha
u'(s,\cd)\|_2 \langle s\rangle^{-4/5+\delta}\, ds,
\end{multline*}
if $\delta>0$.  But if $\delta$ is chosen small enough so that
$4/5-\delta >1/2+1/10$ then we can use Schwarz's inequality along
with \eqref{5.11} to see that the last expression is
$O(\varepsilon^2)$.  We are using here the fact that $u=u_0+w$, as
well as the fact that $u_0$  satisfies better bounds than those in
\eqref{5.11} because of \eqref{5.4}.

Let us see that the contribution of the second term in the right
side of \eqref{5.16} enjoys the same bound.  For a given $J$ we
can use \eqref{4.8} to see that the contribution is dominated by
\begin{multline}\label{5.17}
\int_0^t\langle s\rangle^{-1}\Bigl(\sum_{\substack{|\alpha|+\mu
\le 51 \\ \mu\le 2}}\|L^\mu Z^\alpha u'(s,\cd)\|_2 +
\sum_{\substack{|\alpha|+\mu \le 50 \\ \mu \le 1}}\|\langle
t+r\rangle L^\mu Z^\alpha \square u(s,\cd)\|_2\Bigr)
\\
\times \sum_{\substack{|\alpha|+\mu\le 50 \\ \mu\le 1}}\|\langle
y\rangle^{-1}L^\mu Z^\alpha u'(s,\cd)\|_2 \, ds
\\
+\int_0^t\sum_{R=2^k<\frac{c_0s}2}\sum_{\substack{|\alpha|+\mu \le
50 \\ \mu\le 1}}\Bigl(\|\langle y\rangle^{-1}L^\mu Z^\alpha
u'(s,\cd)\|_{L^2(|y|\approx R)}+\|\langle y\rangle^{-1}L^\mu
Z^\alpha u(s,\cd)\|_{L^2(|y|\approx R)}\Bigr)
\\
\times \sum_{\substack{|\alpha|+\mu \le 50 \\ \mu\le 1}}\|\langle
y\rangle^{-1}L^\mu Z^\alpha u'(s,\cd)\|_{L^2(|y|\approx R)} \, ds,
\end{multline}
with $c_0=\min_I c_I$, and $L^2(|y|\approx R)$ indicating
$L^2$-norms over $\{y\in \mathbb{R}^3\backslash \mathcal{K}: \,
|y|\in [R/4,2R]\}$.  If one uses \eqref{5.4} and \eqref{5.8} to
estimate the first factor in the last term, one concludes that
this term is dominated by
$$\varepsilon \int_0^t (\log(2+s))^2 \langle
s\rangle^{-4/5}\sum_{\substack{|\alpha|+\mu\le 50 \\ \mu\le
1}}\|\langle y\rangle^{-1/2}L^\mu Z^\alpha u'(s,\cd)\|_2 \, ds
=O(\varepsilon^2),$$ using \eqref{5.11} and \eqref{5.4} in the
last step.  For the first term of \eqref{5.17}, we note that by
\eqref{5.8}
\begin{align}\label{5.18}
\langle s+r\rangle \sum_{\substack{|\alpha|+\mu \le 50 \\ \mu\le
1}}|L^\mu Z^\alpha \square u|
&\le C\langle s+r\rangle
\sum_{\substack{|\alpha|+\mu \le 51 \\ \mu\le 1}}|L^\mu Z^\alpha
u'|^2
\\
&\le C\varepsilon \log(2+s) (1+s)^{1/5}
\sum_{\substack{|\alpha|+\mu \le 51 \\ \mu\le 1}}|L^\mu Z^\alpha
u'|,
\notag
\end{align}
assuming, as we may, that $\varepsilon\le 1$.  Thus, by
\eqref{5.10} and \eqref{5.4} the contribution of the first term in
the right side of \eqref{5.17} must be dominated by
$$\int_0^t \langle s\rangle^{-1}\Bigl(\varepsilon \langle
s\rangle^{1/10}+\varepsilon \log(2+s)\langle
s\rangle^{1/10+1/5}\Bigr)\sum_{\substack{|\alpha|+\mu \le 51 \\
\mu\le 2}}\|\langle y\rangle^{-1}L^\mu Z^\alpha u'(s,\cd)\|_2 \,
ds,$$ which is also $O(\varepsilon^2)$ by \eqref{5.11} and
\eqref{5.4}.

This concludes the proof that the null form terms have
$O(\varepsilon^2)$ contributions to $I$.  If we use \eqref{4.7} it
is clear that the multi-speed quadratic terms
$$\sum_{j,k,l=0}^3 B^{IJ,jk}_{K,l}\partial_l
u^K\partial_j\partial_k u^J, \quad J\ne K$$ will have the same
contribution. This completes the proof that $I$ satisfies the
bounds in \eqref{5.12}.

It is also easy to see now that $II$ is $O(\varepsilon^2)$.  If we
use \eqref{5.18}, we see that $II$ is dominated by
$$\varepsilon\int_0^t \langle
s\rangle^{-4/5}\log(2+s)\sum_{\substack{|\alpha|+\mu \le 51 \\
\mu\le 1}}\|L^\mu Z^\alpha u'(s,\cd)\|_{L^2(|y|<4)}\, ds,$$ which
is $O(\varepsilon^2)$ by \eqref{5.11} and \eqref{5.4}.

This completes step i) of the proof, which was to show that
\eqref{5.8}-\eqref{5.11} imply \eqref{5.7}.

To finish the proof of Theorem \ref{theorem1.1} we need to show
how \eqref{5.7} implies \eqref{5.8}-\eqref{5.12}.  In proving the
$L^2$ estimates we shall use the fact that, in the notation of \S
1, $ \square_\gamma u =B(du), $ where the quadratic form $B(du)$
is the semilinear part of the nonlinearity $Q$, and
$$\gamma^{IJ,jk}=\gamma^{IJ,jk}(u')=-\sum_{\substack{0\le l\le 3 \\ 1\le K\le
D}}B^{IJ,jk}_{K,l}\partial_l u^K.$$
Depending on the linear estimates we shall employ, at times we
shall prove certain $L^2$ bounds for $u$ while at other times, we
shall prove them for $w$. Since $u=w+u_0$ and $u_0$ satisfies the
bounds in \eqref{5.4} it will always be the case that bounds for
$w$ will imply those for $u$ and vice versa.  Also note that by
\eqref{5.7}
\begin{equation}\label{5.19}\|\gamma'(s, \cd)\|_\infty \le
\frac{C\varepsilon}{(1+s)}.
\end{equation}

Using these facts we can prove \eqref{5.9}.  Let us first notice
that if we use \eqref{2.5} and \eqref{5.7} then we can estimate
the energy of $\partial^j_tu$ for $j\le M \le 100$.  We shall use
induction on $M$.

We first notice that by \eqref{2.5} and \eqref{5.19} we have
\begin{equation}\label{5.20}
\partial_t E^{1/2}_M(u)(t)\le C\sum_{j\le M}\|\square_\gamma
\partial^j_t u(t,\cd)\|_2 +
\frac{C\varespilon}{1+t}E_M^{1/2}(u)(t).
\end{equation}
Note that for $M=1,2,\dots$ and
\begin{align*}
\sum_{j\le M}|\square_\gamma \partial^j_tu| &\le C\Bigl(\sum_{j\le
M}|\partial^j_tu'|+\sum_{j\le M-1}|\partial^j_t \partial^2
u|\Bigr)\sum_{|\alpha|\le 40}|\partial^\alpha u'|
\\ &\qquad\qquad\qquad\qquad
+C\sum_{|\alpha|\le M-41}|\partial^\alpha u'|\sum_{40< |\alpha|\le
M/2}|\partial^\alpha u'|
\\
&\le \frac{C\varepsilon}{1+t}\Bigl(\sum_{j\le
M}|\partial^j_tu'|+\sum_{j\le M-1}|\partial^j_t \partial^2
u|\Bigr) + C\sum_{|\alpha|\le M-41}|\partial^\alpha
u'|\sum_{|\alpha|\le M/2}|\partial^\alpha u'|,
\end{align*}
since \eqref{5.7} and \eqref{5.4} imply $|\partial^\alpha u'|\le
C\varepsilon/(1+t)$ if $|\alpha|\le 40$.  Also, if we use elliptic
regularity and repeat this argument we get
\begin{align*}
\sum_{j\le M-1}\|\partial^j_t\partial^2 u(t,\cd)\|_2&\le
C\sum_{j\le M}\|\partial^j_t u'(t,\cd)\|_2 + C\sum_{j\le
M-1}\|\partial^j_t \square u(t,\cd)\|_2
\\
&\le C\sum_{j\le M}\|\partial^j_t
u'(t,\cd)\|_2+\frac{C\varepsilon}{1+t}\sum_{j\le
M-1}\|\partial^j_t\partial^2u(t,\cd)\|_2
\\
&\qquad\qquad+C\sum_{|\alpha|\le M-41, |\beta|\le M/2}
\|\partial^\alpha u'(t,\cd)\partial^\beta u'(t,\cd)\|_2
\end{align*}
If $\varepsilon$ is small we can absorb the second to last term
into the left side of the preceding inequality.  Therefore, if we
combine the last two inequalities we conclude that
\begin{multline*}\sum_{j\le M}\|\square_\gamma \partial^j_t u(t,\cd)\|_2\le
\frac{C\varepsilon}{1+t}\sum_{j\le M}\|\partial^j_tu'(t,\cd)\|_2
\\
+C\sum_{|\alpha|\le M-41, |\beta|\le M/2} \|\partial^\alpha
u'(t,\cd)\partial^\beta u'(t,\cd)\|_2.
\end{multline*}
If we combine this with \eqref{5.20} we get that for small
$\varepsilon>0$
\begin{equation}\label{5.21}
\partial_t E^{1/2}_M(u)(t)\le
\frac{C\varepsilon}{1+t}E^{1/2}_M(u)(t)+ C\sum_{|\alpha|\le M-41,
|\beta|\le M/2} \|\partial^\alpha u'(t,\cd)\partial^\beta
u'(t,\cd)\|_2,
\end{equation}
since when $\varepsilon$ is small $\tfrac12 E^{1/2}_M(u)(t)\le
\sum_{j\le M}\|\partial^j_t u'(t,\cd)\|_2\le 2E^{1/2}_M(u)(t)$.

If $M=40$, the last term in \eqref{5.21} drops out and so
$$\partial_t E^{1/2}_{40}(u)(t)\le
\frac{C\varepsilon}{1+t}E^{1/2}_{40}(u)(t).$$ Since
$E^{1/2}_{100}(u)(0)\le C\varepsilon$, an application of
Gronwall's inequality yields
\begin{equation}\label{5.22}
\sum_{j\le 40}\|\partial^j_tu'(t,\cd)\|_2\le 2E_{40}^{1/2}(u)(t)\le
C\varepsilon (1+t)^{C\varepsilon}.\end{equation}
 By elliptic
regularity and \eqref{5.7} this leads to the bounds
$$\sum_{|\alpha|\le 40}\|\partial^\alpha u'(t,\cd)\|_2\le
C\varepsilon (1+t)^{C\varepsilon}.$$

If $M>40$ we have to deal with the last term in \eqref{5.21}.  To
do this we first note that by Lemma 4.5 we have
$$\sum_{|\alpha|\le M-41,
|\beta|\le M/2} \|\partial^\alpha u'(t,\cd)\partial^\beta
u'(t,\cd)\|_2\le C\sum_{|\gamma|\le \max(M-39, 2+M/2)}\|\langle
x\rangle^{-1/2}Z^\gamma u'(t,\cd)\|_2^2,$$ which means that for
$40<M\le 100$, \eqref{5.21} and Gronwall's inequality yield
\begin{equation}\label{5.23}
E^{1/2}_M(u)(t)\le
C(1+t)^{C\varepsilon}\Bigl[\varepsilon + \sum_{|\alpha|\le \max(M-39,
2+M/2)}\|\langle x\rangle^{-1/2}Z^\alpha u'\|_{L^2(S_t)}^2\Bigr],
\end{equation}
if, as before, $S_t=[0,t]\times \mathbb{R}^3\backslash
\mathcal{K}$.

If we use \eqref{5.22} and \eqref{5.23} along with a simple
induction argument we conclude that we would have the desired
bounds
\begin{equation}\label{5.24}
E^{1/2}_{100}(u)(t)\le C\varepsilon(1+t)^{C\varepsilon+\sigma}
\end{equation}
for arbitrarily small $\sigma>0$ if we could prove the following

\begin{lemma}\label{lemma5.2}  Under the above assumptions if
$M\le 100$ and
\begin{multline}\label{5.25}
\sum_{|\alpha|\le M}\|\partial^\alpha
u'(t,\cd)\|_2+\sum_{|\alpha|\le M-3}\|\langle
x\rangle^{-1/2}\partial^\alpha u'\|_{L^2(S_t)} +\sum_{|\alpha|\le
M-4}\|Z^\alpha u'(t,\cd)\|_2
\\
+\sum_{|\alpha|\le M-6}\|\langle x\rangle^{-1/2}Z^\alpha
u'\|_{L^2(S_t)} \le C\varepsilon (1+t)^{C\varepsilon+\sigma},
\end{multline}
with $\sigma>0$, then there is a constant $C'$ so that
\begin{multline}\label{5.26}
\sum_{|\alpha|\le M-2}\|\langle x\rangle^{-1/2}\partial^\alpha
u'\|_{L^2(S_t)}+\sum_{|\alpha|\le M-3}\|Z^\alpha
u'(t,\cd)\|_2
\\
+\sum_{|\alpha|\le M-5}\|\langle x\rangle^{-1/2}Z^\alpha
u'\|_{L^2(S_t)}\le C'\varepsilon (1+t)^{C'\varepsilon+C'\sigma}.
\end{multline}
\end{lemma}

\noindent{\bf Proof of Lemma \ref{lemma5.2}:}  Let us start out by
estimating the first term in the right side of \eqref{5.26}.  By
\eqref{5.4} and  \eqref{2.21} we have
\begin{align}\label{5.27}
(\log(2+t))^{-1/2}&\sum_{|\alpha|\le M-2}\|\langle
x\rangle^{-1/2}\partial^\alpha u'\|_{L^2(S_t)}
\\
&\le C\varespilon + (\log(2+t))^{-1/2}\sum_{|\alpha|\le
M-2}\|\langle x\rangle^{-1/2}\partial^\alpha w'\|_{L^2(S_t)}
\notag
\\
&\le C\varepsilon + C\sum_{|\alpha|\le M-1}\int_0^t
\|\partial^\alpha \square w(s,\cd)\|_2\, ds \notag
+C\sum_{|\alpha|\le M-2}\|\partial^\alpha \square w\|_{L^2(S_t)}.
\notag
\end{align}
Since $\partial^\alpha \square w=\partial^\alpha\square
u-\partial^\alpha\square u_0$, \eqref{5.4} implies that the right
side is
$$\le C\varepsilon+C\sum_{|\alpha|\le
M-1}\int_0^t\|\partial^\alpha \square u(s,\cd)\|_2ds
+C\sum_{|\alpha|\le M-2}\|\partial^\alpha \square u\|_{L^2(S_t)}.$$
If $M\le 40$ we can use \eqref{5.7} and \eqref{5.25} to see that
the last two terms are $\le
C\varepsilon(1+t)^{C\varepsilon+\sigma}$.  If $40<M\le 100$ we can
repeat the proof of \eqref{5.23} to conclude that they are
\begin{align*}
&\le C\varepsilon(1+t)^{2C\varepsilon+2\sigma}+C\sum_{|\alpha|\le
\max(M-39,2+M/2)}\|\langle x\rangle^{-1/2}Z^\alpha
u'\|_{L^2(S_t)}^2\\
&\quad+C\sup_{0\le s\le t}\Bigl(\sum_{|\alpha|\le M-6} \|Z^\alpha
u'(s,\cd)\|_2\Bigr) \sum_{|\alpha|\le
\max(M-39,2+M/2)}\|\langle x\rangle^{-1/2}Z^\alpha u'\|_{L^2(S_t)}
\\
&\le C\varepsilon (1+t)^{2C\varepsilon+2\sigma},
\end{align*}
using the induction hypothesis \eqref{5.25} and the fact that
$\max (M-39,2+M/2)\le M-6$ if $M\ge 40$.  Thus, the left side of
\eqref{5.27} is $\le C\varepsilon(1+t)^{2C\varepsilon +2\sigma}$,
which by \eqref{5.4} means that
$$\sum_{|\alpha|\le M-2}\|\langle x\rangle^{-1/2}\partial^\alpha
u'\|_{L^2(S_t)}\le C\varepsilon
(1+t)^{2C\varepsilon+2\sigma}\log(2+t).$$ Thus, we have the
desired bounds for the first term in the left side of
\eqref{5.26}.

We need to control the second term in the left side of
\eqref{5.26}.  Here we need to use \eqref{2.19}.  In order to do
so, we need to estimate the first term in its right side.  We note
that if $Y_{M-3,0}(t)$ is as in \eqref{2.19}, then
\begin{align*}
&\sum_{|\alpha|\le M-3}\|\square_\gamma Z^\alpha u(t,\cd)\|_2
\\
&\le C\sum_{|\beta|+|\gamma|\le M-3}\|Z^\beta u'(t,\cd) Z^\gamma
u'(t,\cd)\|_2
\\
&\le C\sum_{|\beta|\le M-3, |\gamma|\le 40}\|Z^\beta u'(t,\cd)\|_2
\|Z^\gamma u'(t,\cd)\|_\infty + C\sum_{|\beta|,|\gamma|\le
M-43}\|Z^\beta u'(t,\cd)Z^\gamma u'(t,\cd)\|_2
\\
&\le \frac{C\varepsilon}{1+t}Y^{1/2}_{M-3,0}(t) + C\sum_{|\beta|\le
M-41}\|\langle x\rangle^{-1/2}Z^\beta u'(t,\cd)\|^2_2.
\end{align*}
In the last step, we used \eqref{5.7} and Lemma \ref{lemma4.5}.
By plugging this into \eqref{2.19}, we conclude that
\begin{multline*}
\partial_t Y_{M-3,0}(t) \le \frac{C\varepsilon}{1+t}Y_{M-3,0}(t) +
C\sum_{|\beta|\le M-41}\|\langle x\rangle^{-1/2} Z^\beta
u'(t,\cd)\|_2^2
\\
+C\sum_{|\alpha|\le M-2}\|\langle x\rangle^{-1/2}\partial^\alpha
u'(t,\cd)\|_2^2.
\end{multline*}
Therefore, by Gronwall's inequality, we have
\begin{multline*}
\sum_{|\alpha|\le M-3}\|Z^\alpha u'(t,\cd)\|^2_2 \le CY_{M-3,0}(t)
\\\le C(1+t)^{C\varepsilon}\Bigl( \varepsilon^2 + C\sum_{|\beta|\le
M-41}\|\langle x\rangle^{-1/2} Z^\beta u'(t,\cd)\|_{L^2(S_t)}^2
\\
+C\sum_{|\alpha|\le M-2}\|\langle x\rangle^{-1/2}\partial^\alpha
u'(t,\cd)\|_{L^2(S_t)}^2\Bigr).
\end{multline*}
In the previous step we estimated the last term in the right.
Since the inductive hypothesis handles the second term, we
conclude that the second term in \eqref{5.26} also satisfies the
desired bounds.  Using \eqref{2.22}, this in turn implies that the
third term satisfies the bounds, which completes the proof. \qed


This proves \eqref{5.24}.  By elliptic regularity, we get
$$\sum_{|\alpha|\le 100}\|\partial^\alpha u'(t,\cd)\|_2\le
C\varepsilon(1+t)^{C\varepsilon+\sigma},$$ which in turn yields
\eqref{5.9}.  We also get from Lemma \ref{lemma5.2} that
\begin{multline}\label{5.28}
\sum_{|\alpha|\le 98}\|\langle x\rangle^{-1/2}\partial^\alpha
w'\|_{L^2(S_t)}+\sum_{|\alpha|\le 97}\|Z^\alpha w'(t,\cd)\|_2
\\
+\sum_{|\alpha|\le 95}\|\langle x\rangle^{-1/2}Z^\alpha
w'\|_{L^2(S_t)}\le C'\varepsilon (1+t)^{C'\varepsilon+C'\sigma},
\end{multline}
since the same sort of bounds hold when $u$ is replaced by $w$.

Here and in what follows $\sigma$ denotes a small constant that
must be taken to be larger and larger at each occurrence.  Note
that in terms of the number of $Z$ derivatives \eqref{5.26} is
considerably stronger than the variants of \eqref{5.10} and
\eqref{5.11} where one just takes the terms with $\nu=0$.  This is
because just as in going from \eqref{5.9} to \eqref{5.28} there is
a loss of derivatives, there will be a loss of derivatives in
going from $L^2$ bounds for terms of the form $L^\nu Z^\alpha w'$
to those of the form $L^{\nu+1}Z^\alpha w'$.
%

The proof of the estimates involving powers of $L$ is a bit more
complicated.  Still we shall follow the above strategy.  First we
shall estimate $L^\nu \partial^\alpha u'$ in $L^2$ when $\alpha$
is small using \eqref{5.7}.  Then we shall estimate the remaining
parts of \eqref{5.10} and \eqref{5.11} for this value of $\nu$ by
an inductive argument that is similar to the one in Lemma
\ref{lemma5.2}.

The main part of the next step will be to show that
\begin{equation}\label{5.29} \sum_{\substack{|\alpha|+\mu\le 92 \\
\mu\le 1}}\|L^\mu \partial^\alpha u'(t,\cd)\|_2\le C\varepsilon
(1+t)^{C\varepsilon+\sigma}.
\end{equation}

    For this we shall want to use \eqref{2.16}.  We first must
establish appropriate versions of \eqref{2.15} for $N_0+\nu_0\le
92$, $\nu_0=1$.  For this we note that for $M\le 92$
\begin{align*}
\sum_{\substack{j+\mu\le M \\ \mu\le 1}}\Bigl(|\tilde L^\mu&
\partial_t^j \Box_\gamma u|+|[\tilde L^\mu \partial^j_t,
\square -\square_\gamma]u|\Bigr)
\\
 &\le C\Bigl(\sum_{j\le M-1}|\tilde
L\partial_t^j u'|+\sum_{j\le M-2}|\tilde L \partial_t^j \partial^2
u|\Bigr) \sum_{|\alpha|\le 40}|\partial^\alpha u'|
\\
&+C\sum_{|\alpha|\le M-41}|L\partial^\alpha u'|\sum_{|\alpha|\le
M}|\partial^\alpha u'| +C\sum_{|\alpha|\le M}|\partial^\alpha
u'|\sum_{|\alpha|\le \max(M/2,M-40)}|\partial^\alpha u'|.
\end{align*}
From this, \eqref{5.7}, Lemma 4.5 and elliptic regularity we get
that for $M\le 92$
\begin{align*}\sum_{\substack{j+\mu \le M\\ \mu\le
1}}\Bigl(\|\tilde L^\mu \partial_t^j \square_\gamma u(t,\cd)&\|_2
+\|[\tilde L^\mu \partial^j_t, \square
-\square_\gamma]u(t,\cd)\|_2\Bigr)
\\
&\le \frac{C\varepsilon}{1+t}\sum_{\substack{j+\mu \le M\\ \mu\le
1}}\|\tilde L^\mu \partial^j_t u'(t,\cd)\|_2
\\
&+C\sum_{|\alpha|\le M-41}\|\langle
x\rangle^{-1/2}L\partial^\alpha u'(t,\cd)\|_2 \sum_{|\alpha|\le
94}\|\langle x\rangle^{-1/2}Z^\alpha u'(t,\cd)\|_2
\\
&+C\sum_{|\alpha|\le \max (M, 2+M/2)}\|\langle
x\rangle^{-1/2}Z^\alpha u'(t,\cd)\|_2^2.
\end{align*}
Based on this if $\varepsilon$ is small then \eqref{2.15} holds
with $\delta =C\varepsilon$ and
$$H_{1,M-1}(t)=C\sum_{|\alpha|\le M-41}\|\langle x\rangle^{-1/2}L
\partial^\alpha u'(t,\cd)\|_2^2 + C\sum_{|\alpha|\le 94}\|\langle
x\rangle^{-1/2}Z^\alpha u'(t,\cd)\|_2^2.$$ Therefore since the
conditions on the data give $X_{\mu,j}(0)\le C\varepsilon$ if
$\mu+j\le 100$ it follows from \eqref{2.16} and \eqref{5.28} that
for $M\le 92$
\begin{align}\label{5.30}
\sum_{\substack{|\alpha|+\mu\le M \\ \mu\le 1}}\|L^\mu
\partial^\alpha u'(t,\cd)\|_2&\le
C\varepsilon(1+t)^{C\varepsilon+\sigma}+C(1+t)^{C\varepsilon}\sum_{|\alpha|\le
M-41}\|\langle x\rangle^{-1/2}L\partial^\alpha u'\|_{L^2(S_t)}^2
\\
&+C(1+t)^{C\varepsilon}\int_0^t\sum_{|\alpha|\le
M+1}\|\partial^\alpha u'(s,\cd)\|_{L^2(|x|<1)} \, ds. \notag
\end{align}
If we apply \eqref{2.32} and \eqref{5.4} we get that the last
integral is dominated by $\varepsilon\log(2+t)$ plus
\begin{multline*}\int_0^t \sum_{|\alpha|\le M+1}\|\partial^\alpha
w'(s,\cd)\|_{L^2(|x|<1)}ds
\\
\le C\sum_{|\alpha|\le M+2}\int_0^t\Bigl(\int_0^s\|\partial^\alpha
\square w(\tau,\cd)\|_{L^2(|\, |x|-(s-\tau) \,
|<10)}d\tau\Bigr)ds.\end{multline*} By \eqref{5.4} if we replace
$w$ by $u_0$ then the analog of the last term is
$O(\log(2+t)\varepsilon)$.  We therefore conclude that
\begin{multline*}\sum_{|\alpha|\le M+1}\int_0^t\|\partial^\alpha
u'(s,\cd)\|_{L^2(|x|<1)}ds \le C\log(2+t)\varepsilon
\\
+C\sum_{|\alpha|\le M+2}\int_0^t\Bigl(\int_0^s\|\partial^\alpha
\square u(\tau,\cd)\|_{L^2(|\, |x|-(s-\tau)\,
|<10)}d\tau\Bigr)ds.\end{multline*}
 Since
 $$\sum_{|\alpha|\le M+2}|\partial^\alpha \square u|\le
 C\sum_{|\alpha|\le M+3}|\partial^\alpha u'|\sum_{|\alpha|\le
 1+M/2}|\partial^\alpha u'|,$$
 an application of Lemma \ref{lemma4.5} yields
 $$\sum_{|\alpha|\le M+2}\|\partial^\alpha \square
 u(\tau,\cd)\|_{L^2(|\, |x|-(s-\tau)\, |<10)}\le
 C\sum_{|\alpha|\le 95}\|\langle x\rangle^{-1/2}Z^\alpha
 u'\|_{L^2(|\, |x|-(s-\tau)\, |<20)}^2,$$
 since $3+M/2\le 95$ if $M\le 92$.  Since the sets $\{(\tau,x): \,
 |\, |x|-(j-\tau)\, |<20\}$, $j=0,1,2,\dots$ have finite overlap,
 we conclude that for $M\le 92$
\begin{align*}
\sum_{|\alpha|\le M+1}\int_0^t\|\partial^\alpha
 u'(s,\cd)\|_{L^2(|x|<1)}ds &\le C\varepsilon \log(2+t)+C\sum_{|\alpha|\le
95}\|\langle
 x\rangle^{-1/2}Z^\alpha u'\|_{L^2(S_t)}^2
\\
 &\le
 C\varespilon(1+t)^{C\varepsilon+\sigma}.\end{align*}
 Therefore, by \eqref{5.30} we have that
 \begin{align*}
 \sum_{\substack{|\alpha|+\mu \le M \\ \mu\le 1}}\|L^\mu
 \partial^\alpha u'(t,\cd)\|_2 &\le
 C\varepsilon(1+t)^{C\varepsilon+\sigma}
\\
 &+C(1+t)^{C\varepsilon}\sum_{|\alpha|\le
 M-41}\|\langle x\rangle^{-1/2}L \partial^\alpha
 u'\|_{L^2(S_t)}^2.\end{align*}
 This gives the desired bounds when $M\le 40$.

 If we now use \eqref{2.21} with $\nu_0=1$ and $N_0+\nu_0=92$,
 then the analog of Lemma \ref{lemma5.2} where $M=100$ is replaced
 by $M=92$ and $u$ is replaced by $Lu$ is valid.  By an induction
 argument we get \eqref{5.29} from this as well as
 \begin{multline}\label{5.31}
 \sum_{\substack{|\alpha|+\mu\le 90\\ \mu\le 1}}\|\langle
 x\rangle^{-1/2}L^\mu \partial^\alpha
 w'\|_{L^2(S_t)}+\sum_{\substack{|\alpha|+\mu\le 89 \\ \mu\le
 1}}\|L^\mu Z^\alpha w'(t,\cd)\|_2
 \\
 +\sum_{\substack{|\alpha|+\mu\le 87 \\ \mu\le 1}}\|\langle
 x\rangle^{-1/2}L^\mu Z^\alpha w'\|_{L^2(S_t)}\le
 C\varepsilon(1+t)^{C\varepsilon+C\sigma}.\end{multline}

 If we repeat this argument we can estimate $L^2Z^\alpha u'$ and
 $L^3Z^\alpha u'$ for appropriate $Z^\alpha$.  Using \eqref{5.29}
 and \eqref{5.31} and the last argument gives
 \begin{multline*}\sum_{\substack{|\alpha|+\mu\le 84 \\ \mu\le
 2}}\|L^\mu \partial^\alpha
 w'(t,\cd)\|_2+\sum_{\substack{|\alpha|+\mu \le 81 \\ \mu\le
 2}}\|L^\mu Z^\alpha w'(t,\cd)\|_2
 \\
 +\sum_{\substack{|\alpha|+\mu \le 79 \\ \mu\le 2}}\|\langle
 x\rangle^{-1/2}L^\mu Z^\alpha w'\|_{L^2(S_t)}\le
 C\varepsilon(1+t)^{C\varepsilon+C\sigma}.
 \end{multline*}
 Then using the estimates for $L^\mu Z^\alpha u'$, $\mu\le 2$ we
 can argue as above to finally get
 \begin{multline*}
 \sum_{\substack{|\alpha|+\mu \le 76 \\ \mu\le 3}}\|L^\mu
 \partial^\alpha w'(t,\cd)\|_2 +\sum_{\substack{|\alpha|+\mu \le
 73 \\ \mu\le 3}}\|L^\mu Z^\alpha w'(t,\cd)\|_2
 \\
 +\sum_{\substack{|\alpha|+\mu\le 70 \\ \mu\le 3}}\|\langle
 x\rangle^{-1/2}L^\mu Z^\alpha w'\|_{L^2(S_t)}\le
 C\varepsilon(1+t)^{C\varepsilon+C\sigma}.
 \end{multline*}
 If we combine this with our earlier bounds, we conclude that
 \eqref{5.10} and \eqref{5.11} must be valid.

 It remains to prove \eqref{5.8}.  This is straightforward.  If we
 use Theorem \ref{theorem3.1} we find that its left side is
 dominated by the square of that of \eqref{5.11}.  Hence
 \eqref{5.11} implies \eqref{5.8}, which finishes the proof.

\end{document}